\documentclass[12pt]{article}
\usepackage{graphicx}
\usepackage{amssymb,amsmath}
\usepackage{hyperref}
\usepackage[margin=1in]{geometry}
\usepackage{mathrsfs} 
\usepackage{tikz}
\usetikzlibrary{decorations.pathmorphing}
\usetikzlibrary{arrows.meta, positioning, fit}

\usepackage{float}          
\usepackage{algorithm}      
\usepackage{algpseudocode}  

\definecolor{blue1}{RGB}{173,216,240} 
\definecolor{blue2}{RGB}{186,234,250} 
\newcommand{\defeq}{\mathrel{\mathop:}=}
\newcommand{\R}{\mathbb{R}}

\newcommand{\rev}[1]{\textcolor{red}{#1}}
\newcommand{\iu}{\mathrm{i}\mkern1mu}

\title{A Neural Network Enhanced Born Approximation for Inverse Scattering
}
\author{
Ansh Desai\footnotemark[1]
\and Jonathan Ma\thanks{Mathematical Sciences, University of Delaware, Newark DE 19716.
email:{\tt adesai@udel.edu}, email:{ \tt johnma@udel.edu}}
\and Timo L\"ahivaara\thanks{Department of Technical Physics, University of Eastern Finland, Kuopio, Finland. 
email:{\tt timo.lahivaara@uef.fi}} 
\and Peter Monk\thanks{Mathematical Sciences, University of Delaware, Newark DE 19716. email: {\tt monk@udel.edu}}
}

\date{}

\begin{document}

\maketitle
\makeatletter
\renewcommand{\@makefnmark}{} 
\footnotetext{%
\noindent {\bf{}Funding}: The research of A.D. is supported by the University of Delaware Undergraduate Research Program.
The research of T.L. is partially supported by the Research Council of Finland via the Finnish Center of Excellence of Inverse Modeling and Imaging, 
the research project 321761, and the Flagship of Advanced Mathematics for Sensing Imaging and Modeling grant 358944.
The research of P.M. was partially supported by the US AFOSR under grant number FA9550-23-1-0256.
}
\makeatother

\begin{abstract}
    Time-harmonic acoustic inverse scattering concerns the ill-posed and nonlinear problem of determining the refractive index of an inaccessible, penetrable scatterer based on far field wave scattering data. 
    When the scattering is weak, the regularized inverse Born approximation provides a linearized model for recovering the shape and material properties of a scatterer. We propose two convolutional neural network (CNN) algorithms to correct the traditional inverse Born approximation even when the scattering is not weak. These are denoted Born-CNN (BCNN) and CNN-Born (CNNB). BCNN applies a post-correction to the Born 
    reconstruction, while CNNB pre-corrects the data. Both methods leverage the Born approximation's excellent fidelity in weak scattering, while extending its applicability beyond its theoretical limits. CNNB particularly exhibits a strong generalization to more complex out of distribution scatterers.
    Based on numerical tests and benchmarking against other standard approaches, our corrected Born models provide alternative data-driven 
    methods 
    for obtaining the refractive index, extending the utility of the Born approximation to regimes where the traditional method fails. 
\end{abstract}



\section{Introduction}
The time-harmonic inverse scattering problem concerns the determination of the properties, such as shape or material composition, of an inaccessible object from remote measurements of wave scattering data.  Problems of this type arise in, for example, seismology, medical imaging, and radar applications. Due to the significance of the applications, there has been extensive work on various algorithms for this problem (see for example~\cite{colton+kress4}).  

We shall study a particular inverse scattering problem: the inverse medium problem for the Helmholtz equation.  In this case, we seek to reconstruct the refractive index 
of a bounded penetrable scatterer from far field acoustic data.  
The major complication is that this inverse scattering problem is both ill-posed and nonlinear~\cite[Theorem 4.21, page 448]{colton+kress4}. Current approaches can be divided into two broad classes: 1) quantitative methods and 2) qualitative methods.
A quantitative method attempts to reconstruct the scatterer
directly. Often this involves a non-linear and regularized optimization problem that seeks a reconstruction
corresponding to a far field pattern that matches measurement data.  This is computationally intensive and can fail due to local minima~\cite{colton+kress4}.  Another quantitative approach is to assume that the scattering is weak and to use the inverse Born approximation.  We shall discuss this in more detail shortly since correcting this approach is the subject of our paper. In contrast to these quantitative methods, the goal of a qualitative method such as the Linear Sampling
Method~\cite{CCH,colton+kress4} is to approximate the support or boundary of the scatterer.  Approaches of this type do not involve optimization or the solution of the forward problem.  They do not require strong a priori assumptions, but cannot directly distinguish different materials in the field of view.  

Turning now to the Born approximation, the forward scattering problem that maps a known refractive index to the predicted far field pattern is given by a Neumann series in a certain integral operator, provided an appropriate norm of the integral operator is less than one.  Conditions for this to occur have been derived in several cases (for example,
\cite{Hudson:81} in the seismic context) and we recall a simple sufficient condition in
Section~\ref{sec:iscat}. Selecting only the first term in this series defines a linear map, the Born approximation, from the contrast to the approximate far field pattern.

For the inverse problem, the inverse Born approximation can be understood as inverting the
linearized Born approximation. Inverting this linear map is ill-posed, but removes the difficulty of nonlinearity and provides an avenue to solve the inverse problem approximately using regularization to restore stability.  An early reference for this technique is ~\cite{Born_Wolf:59}, while  applications and computational techniques are described in \cite{Devaney:82,Devaney:12}. The inverse Born approximation has been widely applied to neutron scattering, medical imaging, and seismic inverse problems~\cite{Hudson:81}. Within the weak scattering approximation, there has been a great deal of work to incorporate higher-order expansions into the forward and inverse Born approximation (see the review of Moskow and Schotland~\cite{Moskow:19}).  Our goal in this paper is to
correct the Born approximation and extend its applicability beyond the weak scattering limit using neural networks (NNs) as correctors.

While traditional inverse scattering techniques offer interpretability and theoretical grounding, they often suffer from high computational cost or restrictive assumptions (e.g., weak scattering). To address these limitations, recent years have seen rapid growth in the use of machine learning, particularly neural networks, to learn the inverse map directly or to enhance physics-based approximations.

One machine learning approach to correcting the Born approximation is the statistical approximation error correction used in~\cite{Kaipio:19}, which shows that the Born approximation can be extended outside the weak scattering approximation by a suitable training approach if the scatterer is well represented in the training data.  The approach of~\cite{Kaipio:19} is based on Bayesian statistics and not on neural networks.

We now turn to a review of recent efforts to use NN based  machine learning techniques to address the inverse scattering problem. These approaches vary in architecture and in whether they incorporate aspects of the underlying physics of wave propagation. We emphasize that this is not intended to be a comprehensive review of this rapidly evolving field, but rather a focused overview of representative approaches most relevant to our work.

Viewing the map from the far field to the regularized reconstruction of the scatterer as an unknown nonlinear function, we can appeal to the 
universal approximation property (first proved in \cite{cybenko89}) of NNs to approximate this map. Note that regularization is required because, without it, the map is not defined~\cite{colton+kress4}.

For the inverse problem, we can distinguish two general approaches: 1) Use a NN to approximate the forward problem and then use this model in a traditional optimization framework for the inverse problem (see for example \cite{Stanziola_2021,Saba,Chen:24,liu2024neumannseriesbasedneuraloperator}), 2) Develop end-to-end algorithms that seek to solve the inverse problem directly by approximating the above mentioned data to coefficient map.

We are interested in the second approach.  Within the end-to-end paradigm, an important approach is based on Physics Informed Neural Networks (PINNs)~\cite{Saba,yuyao_2019,yuyao_2022, RAISSI2019686}.  Here the governing differential equation is incorporated in the training. Although it is possible to solve
infinite domain scattering problems using PINNs (for example \cite{wu2022}), we are not aware of a PINNs approach to solving the inverse scattering problem using far field data. Instead, we shall use the
inverse Born approximation to supply  physics information to the networks we design.

Apart from PINNs, there is an extensive literature on end-to-end inversion using supervised NNs including \cite{lexing2019,gao21,Fan:22,Bazow,NIO,Zhou:24,ZHANG24}. In particular, this is the approach in \cite{Fan:22} which is inspired by the Born approximation and is a motivating paper for our work.  They use their BCR-Net~\cite{FAN20191} to approximate the inverse problem since it is based on a nonstandard form of a wavelet representation described in \cite{BCR} and designed to approximate nonlinear integral operators.
We
shall compare our models to the BCR-Net approach in Section~\ref{sec:Num}.

Another end-to-end network motivated by the regularized inverse Born approximation for weak scattering data is considered in \cite{Zhou:24}. That paper also features a mathematical analysis of the generalization and approximation error, and numerical results showing good approximation on weak scattering data.

More generally, beyond architectures motivated by a specific physics-based representation, there has been increasing interest in neural operator frameworks such as the Fourier Neural Operator (FNO) and Deep Operator Network (DeepONet). These
approaches learn mappings between infinite-dimensional function spaces, enabling mesh-independent generalization and rapid solution of PDE-based forward and inverse problems \cite{DBLP:journals/corr/abs-2108-08481, li2020fourier, li2023physicsinformedneuraloperatorlearning, Lu_2021, wang2021learningsolutionoperatorparametric}.  Building on these ideas, a recent hybrid framework termed the Neural Inverse Operator (NIO) combining FNOs and DeepONets was introduced by Molinaro et al.~\cite{NIO}. This has been applied to a wide range of inverse coefficient identification problems, showing excellent results. However, the inverse scattering problem considered in that paper uses a bounded domain and Dirichlet and Neumann data pairs, so it is not directly comparable to our upcoming model using far field data. Nonetheless, with small changes, we can train NIO on far field data and compare results with our Born based schemes. Note that these methods may suffer from hallucinations~\cite{siam-rev}.

Our approach is inspired by the strengths and limitations of the Born approximation, which provides valuable physics-informed structure, especially in the weak scattering regime, but fails for stronger contrasts. Rather than discarding it entirely, we propose to integrate it into supervised learning pipelines as an inductive bias, thereby combining the interpretability of physics-based models with the flexibility of neural networks.  It is also hoped that this approach will
reduce the tendency to hallucinate noted in \cite{siam-rev}.

Our goal is to broaden the applicability of
the Born approximation and arrive at a direct quantitative estimate of the scatterer without the need to compute the solution of a nonlinear inverse problem. We investigate two approaches using a supervised convolutional neural network (CNN) scheme to correct the Born approximation. In particular, we use a generic CNN, but optimize it for each case studied. The first approach--termed \textbf{Born-CNN (BCNN)}--performs a regularized inversion of the Born approximation applied to the far field scattering data and then uses a trained CNN to correct the resulting image. The second approach--termed \textbf{CNN-Born (CNNB)}--trains the neural network to pre-correct the scattering data and then applies a regularized inversion of the Born approximation to produce a corrected reconstruction.  

The novelty of this paper is to suggest and  test the above mentioned schemes combining the Born approximation with a CNN to extend the domain of applicability of the
Born approximation.  In particular, we compare the resulting predictions to those from the standard regularized inverse Born approximation, to a slightly modified version of the Neural Inverse Operator (NIO) from \cite{NIO},
to the BCR-Net network described in \cite{FAN20191}, and to a basic CNN.  After optimizing the networks and training on simple data generated by a few circular scatterers, we
test the networks on several out of distribution scenarios.

The study suggests that combining CNNs and the Born approximation has promise in solving inverse scattering problems and that the training phase does not need to include close copies of the scatterers.

The remainder of this paper proceeds as follows. In Section~\ref{sec:iscat}, we briefly outline the forward and inverse problems underlying this study, summarize the Born approximation, and detail our discretization. Then, in Section~\ref{sec:NN}, we give details of the architectures for the three models considered in the paper and discuss training and testing.  Data generation and the main numerical results of the paper are given in Section \ref{sec:Num}. We draw some conclusions in Section~\ref{sec:concl}.

\section{The Inverse scattering problem}\label{sec:iscat}
In this section we summarize the forward and inverse problems considered in this paper, 
state the Born approximation, and derive our problem setup. The model problem we shall consider is time-harmonic scattering from a penetrable medium modeled by the Helmholtz equation in $\mathbb{R}^2$.
In the upcoming discussion, $k$ denotes the wave number of the field in free space, and $\iu=\sqrt{-1}$.

We suppose that a known incident plane wave with angle of propagation $\varphi$ is given by
\begin{equation}
    u^i(\mathbf{x};\mathbf{d})=\exp(\iu k\mathbf{x}\cdot \mathbf{d}),
\quad \mathbf{d}=\langle \cos\varphi,\sin\varphi\rangle
\end{equation} and that it strikes a bounded penetrable scatterer. 

The square of the refractive index for the medium in which the wave propagates is denoted $\eta(\mathbf{x})\in L^{\infty}(\mathbb{R}^2)$.  This bounded function is assumed to satisfy $\Re (\eta(\mathbf{x}))\geq \eta_{\rm min}>0$, where $\eta_{\rm min}$ is a constant, and $\Im(\eta(\mathbf{x}))\geq 0$ for almost all $\textbf{x}\in\mathbb{R}^2$. In addition, we assume that $|\eta(\mathbf{x})|\leq \eta_{\rm max} $ a.e. for $\mathbf{x}\in\mathbb{R}^2$
 where $\eta_{\rm max}$ is a constant.  The boundedness of the scatterer implies that the \emph{contrast} $\mu(\textbf{x})\defeq \eta(\textbf{x})-1=0$ if $|\textbf{x}|>R$ for some $R>0$ (see, for example \cite{colton+kress4,Kirsch}). 
 
 For a given wave number $k$, the total field $u\defeq u(\mathbf{x};\mathbf{d})$ and the scattered field $u^s\defeq u^s(\mathbf{x};\mathbf{d})$ satisfy the \textbf{Helmholtz equation}:
\begin{align}
    \label{Helmholtz}
    \Delta u+k^2\eta(\mathbf{x})u=0 \text{ in }&\mathbb{R}^2, \\
    u=u^i+u^s \text{ in }&\mathbb{R}^2,
\end{align}
together with the \textbf{Sommerfeld radiation condition}
\begin{align}
    \label{sommerfeld}
    r^{1/2}\left(\frac{\partial u^s}{\partial r}-\iu ku^s\right)&\to 0 \text{ as }r\defeq |\mathbf{x}|\to \infty,
\end{align}
uniformly in $\hat{\mathbf{x}}:=\mathbf{x}/| \mathbf{x}|\in \mathbb{S}^1$ where $\mathbb{S}^1$ denotes the unit circle in $\R^2$.
\begin{figure}
\centering
\begin{tikzpicture}
    \draw[dashed] (-2.5,-2.5) rectangle (2.5,2.5);
    \filldraw[fill=gray!5, draw=none] (-2.5,-2.5) rectangle (2.5,2.5);

    \fill[blue1] (-0.8,0.8) circle [radius=0.7];
    \fill[blue2] (0.8,-0.8) circle [radius=0.5];

    \draw[->,decorate,decoration={snake,amplitude=.5mm,segment length=2mm}] (-3.5,0) -- (-2,0) node[pos=0.4, above] {$u^i$};

    \draw[->,decorate,decoration={snake,amplitude=.5mm,segment length=2mm}] (2,2) -- (3,3) node[midway, above] {$u^s$};

    \node at (1.7,-1.7) {$\Omega$};
    \node at (0,0) {$\eta(\textbf{x})$};
\end{tikzpicture}
\caption{An illustration of the scattering problem. The incident field $u^i$ travels with direction of propagation $\mathbf{d}$ and interacts with the compactly supported scatterer given by $\eta(\textbf{x})\not=1$ shown in blue, and contained in a bounded set $\Omega$ (shown in gray). The scattered field $u^s$ propagates outward with direction angle  $\theta$.}
\label{fig:cartoon}
\end{figure}

Under the conditions given above, Equations (\ref{Helmholtz})-(\ref{sommerfeld}) have a
unique solution for any $k>0$~\cite{colton+kress4}.  The \emph{Forward Problem} consists of solving the above linear well-posed problem given $\mathbf{d}$, $k$, and $\eta$ (see Fig.~\ref{fig:cartoon}).
It follows from the fact that $u^s$ satisfies the Helmholtz equation and the radiation condition that $u^s$ exhibits an asymptotic expansion as an outgoing cylindrical wave for $|\mathbf{x}|$ sufficiently large: 
\begin{equation}
\label{cylindrical}
u^s(\textbf{x};\textbf{d})=\frac{\exp(\iu k|\textbf{x}|)}{\sqrt{|\textbf{x}|}}\left(u_{\infty}(\hat{\textbf{x}},\textbf{d})+\mathcal{O}(|\textbf{x}|^{-1})\right) \text{ as }|\textbf{x}|\to\infty,
\end{equation}
where $u_{\infty}:\mathbb{S}^1\times \mathbb{S}^1\to \mathbb{C}$ is called the \textit{far field pattern} of the scattered wave~\cite{colton+kress4}.

The \emph{Inverse Problem} that we wish to solve is to determine $\mu(\textbf{x})$ (equivalently, $\eta(\mathbf{x})$, the square of the refractive index) given
the far field pattern $u_{\infty}(\hat{\mathbf{x}},\mathbf{d})$ for all $\hat{\mathbf{x}}$ and \textbf{d}
on $\mathbb{S}^1$ (in practice, only a finite set of $\hat{\mathbf{x}}$ and \textbf{d} are used).  Here we assume data is given for a single fixed wave number $k>0$, and that the support of $\mu$ is a priori known to lie in a  bounded search region $\Omega$  (see Fig.~\ref{fig:cartoon}).  This problem is non-linear and ill-posed~\cite{colton+kress4}.

In this study, we use synthetic scattering data generated through a standard finite element approach to approximating (\ref{Helmholtz})-(\ref{sommerfeld}). For each incident direction $\textbf{d}$, the total field $u$ in a neighborhood of the scatterer is computed using the Netgen package~\cite{netgen} with 4th-order finite elements and a mesh-size request of one-eighth of the local wavelength of the wave. The boundary of each scatterer is fitted using isoparametric curved elements. The Sommerfeld radiation condition is handled through a radial Perfectly Matched Layer (PML) implemented by Netgen using complex stretching as discussed in \cite{chew94a}. To generate an approximate far field pattern, we follow \cite{Monk+Suli:98} to map the near field to a far field pattern. The implementation of our forward solver can be found at https://github.com/nibj/Helmholtz-scattering-data.

\subsection{The Born approximation}

Let $H^1_{loc}(\R^2)$ denote the local Sobolev space defined by 
\begin{equation}
\label{sobolev}
    H^{1}_{loc}(\R^2)=\left\{u:\R^2\to \mathbb{C}\:|\:u|_{K}\in H^1(K) \:\forall \mbox{ compact subsets }K\subset\mathbb{R}^2\right\}.
\end{equation}

In particular, let $\Omega$ denote a compact domain 
containing the support of $\mu$ (the sets where $\mu\not=0$). 
It can be shown that if $u\in H^1_{loc}(\R^2)$ is a solution to the scattering problem (\ref{Helmholtz})-(\ref{sommerfeld}), then $u|_{\Omega}\in L^2(\Omega)$ and satisfies the Lippmann-Schwinger equation~\cite{colton+kress4}:
\begin{equation}
    \label{lip-schw}
u(\textbf{x};\textbf{d})=u^i(\textbf{x};\textbf{d})+k^2\int_{\Omega}\Phi(\textbf{x},\textbf{y})\mu(\textbf{y})u(\textbf{y};\textbf{d})\,d\textbf{y} \defeq u^i(\textbf{x};\textbf{d})-(Tu)(\textbf{x};\textbf{d})\quad \forall\:\textbf{x}\in \R^2,\forall\:\textbf{d}\in\mathbb{S}^1
\end{equation}
where $T:L^2(\Omega)\to L^2(\Omega)$ and $$\Phi(\textbf{x},\textbf{y})=\frac{\iu}{4}H^{(1)}_0(k|\mathbf{x}-\mathbf{y}|)$$ is the  fundamental solution to the Helmholtz equation  (\ref{Helmholtz}) and (\ref{sommerfeld}).  Here $H_0^{(1)}(\cdot)$ is the Hankel function of the first kind and order zero. Upon inverting (\ref{lip-schw}), it follows that $u(\textbf{x})$ has a Neumann series representation
\begin{equation}
    \label{neumann}
    u(\textbf{x};\textbf{d})=\sum_{j=0}^{\infty}(-1)^{j}(T^ju^i)(\textbf{x};\textbf{d})
\end{equation}
provided that the operator norm $\|T\|_{2}<1$. A sufficient condition for this is~\cite{schotland+moskow}
\begin{equation}
\label{neumann_criterion}
    \delta\|\mu\|_{\infty}<1,
\end{equation}
where $\|\cdot\|_{\infty}$ is the max norm  and $B_R$ is the smallest ball, with  radius $R$, containing the support of the scatterers and
\[
\delta=k^2\sup_{\textbf{x}\in B_R}\int_{B_R}|\Phi(\textbf{x},\textbf{y})|\,d\textbf{y}.
\]
In $\mathbb{R}^2$, a closed form for $\delta$ as a function of $R$ seems difficult to obtain, but for $k=16$ and $R=\sqrt{2}$ as used in our upcoming numerical results, we can compute numerically that $\| \mu\|_{\infty}<0.00004343641649$.  As we shall see this inequality is is far from necessary for the scatterers we shall use.

When $\|T\|_{2}<1$, the first two terms of (\ref{neumann}) provide the \textbf{Born approximation of the field }$u$:
\begin{equation}
    \label{bornu}
    u(\textbf{x};\textbf{d})\approx u^i(\textbf{x};\textbf{d})+k^2\int_{\Omega}\Phi(\textbf{x},\textbf{y})\mu(\textbf{y})u^i(\textbf{y};\textbf{d})\,d\textbf{y},\quad \forall\:\textbf{x}\in \R^2, \forall\:\textbf{d}\in \mathbb{S}^1.
\end{equation}
Consequently, the far field pattern can be approximated through an asymptotic analysis of (\ref{bornu}) to obtain the \textbf{Born approximation of the far field pattern}:
\begin{equation}
    \label{born}
    u_{\infty}(\hat{\textbf{x}},\textbf{d})\approx \sqrt{\frac{k^3}{8\pi }}\exp\left(\frac{\iu\pi}{4}\right)\int_{\Omega}\exp(\iu k(\textbf{d}-\hat{\textbf{x}})\cdot \textbf{y})\mu(\textbf{y})\,d\textbf{y}\defeq (\mathscr{B}\mu)(\hat{\textbf{x}},\textbf{d}).
\end{equation}
Because of the Neumann series convergence criterion (\ref{neumann_criterion}), the Born approximation is valid for low wave numbers or small contrast, which is termed weak scattering. In particular, the precision of the Born approximation increases as $\|\mu\|_{\infty}$ decreases.

Note that the right hand side of (\ref{born}) is a band-limited Fourier transform of $\eta$, so inverting to find $\eta$ is ill-posed and a regularization technique needs to be used.
A common technique uses Tikhonov regularization and computes the regularized Born approximation
of $\mu$ by
\begin{equation}
    \mu_{\gamma}=({\mathscr B}^*{\mathscr B}+\gamma I)^{-1}{\mathscr B}^*u_{\infty},\label{Bgamma}
\end{equation}
where ${\mathscr B}^*$ is the $L^2$-adjoint of ${\mathscr B}$, $I$ is the identity operator and $\gamma>0$ is a fixed 
regularization parameter determined a priori. Note that
fast methods exist to compute this approximation (including the NN approach of Zhou \cite{Zhou:24} or the low rank approximation method of \cite{Meng25}) but we do not use them here. 

\section{Discretization}
In this paper, we modify the above approach by giving some details of inverting the \textit{stacked Born operator} $\mathcal{B}$ defined in Section~\ref{sec:prob}.
\subsection{Discretization of the Born approximation}\label{sec:prob}
Through translation and rescaling, we assume a priori that the unknown scatterers lie in the square domain $\Omega=[-1,1]^2$. For numerical approximation, this domain is subdivided into an
$N_g\times N_g$ uniform grid of nodal values with coordinates 
$$y^{(1)}_p=-1+\frac{2p}{N_g-1},\quad p=0,1,\ldots,N_g-1,$$
$$y^{(2)}_q=-1+\frac{2q}{N_g-1},\quad q=0,1,\ldots,N_g-1,$$
with $\textbf{y}_{p,q}=(y^{(1)}_p,y^{(2)}_q)$. By a projection onto the grid, we approximate $\mu(\textbf{y})\approx [M(\textbf{y}_{p,q})]\in \R^{N_g\times N_g}$.

We assume a standard source-receiver setup. The incident and scattered fields are  uniformly in $\mathbb{S}^1$, using angles
$$\varphi_i=\frac{2\pi i}{N_S-1}, \quad i=0,1,\ldots,N_S-1,$$
$$\theta_j=\frac{2\pi j}{N_R-1}, \quad j=0,1\ldots,N_R-1,$$
where $N_S$ and $N_R$ denote the number of sources and receivers, respectively. For simplicity, we consider the case in which $N_S=N_R$, following \cite{Fan:22}. Let $\hat{\textbf{x}}_i=\langle \cos\theta_i,\sin\theta_i\rangle$ and $\textbf{d}_j=\langle \cos\varphi_j,\sin\varphi_j\rangle$. By varying both the incident and scattered field directions, it follows that we have available the far field matrix $U_{\infty}\in \mathbb{C}^{N_R\times N_S}$ possibly corrupted by measurement noise or forward modeling errors.

Using quadrature on (\ref{born}), we obtain the components for the discrete Born approximation for a single incident field
\begin{equation}
    \label{Dborn}
    (U_{\infty})_{i,j}\approx \frac{k^2}{\sqrt{8\pi k}}\exp\left(\frac{\iu\pi}{4}\right)h^2\sum_{p,q=0}^{N_g}\exp(\iu k(\textbf{d}_j-\hat{\textbf{x}}_i)\cdot \textbf{y}_{p,q})M_{p,q},
\end{equation}
where $h=2/N_g$. The above approximation in tensor form can be written as a product between the Born $4$-tensor $\mathscr{B}_{\rm disc}\in \mathbb{C}^{N_R\times N_S\times N_g\times N_g}$ and the scatterer matrix $M$
\begin{equation}
\label{Tborn}
    U_{\infty}=\mathscr{B}_{\rm disc}M-\tau_\mu,
\end{equation}
where $\mathscr{B}_{\rm disc}$ encodes the (2+2)-discretization in both angular dimensions and both spatial dimensions, and $\tau_\mu$ is the unknown error due to the Born approximation. 

To utilize linear algebraic methods, we can use standard tensor unfolding methods to collapse $\mathscr{B}_{\rm disc}$ along the spatial and angular dimensions to produce the aforementioned \textit{stacked Born operator} $\mathcal{B}\in \mathbb{C}^{N_R N_S\times N_g^2}$. We thus arrive at the following formulation:
\begin{equation}
\label{stackedBorn}
    \vec{u}_{\infty}=\mathcal{B}\vec{\mu}-\vec{\tau}_\mu,
\end{equation}
where $\vec{u}_{\infty}=\operatorname{vec}(U_{\infty})$, $\vec{\mu}=\operatorname{vec}(M)$, 
and $\vec{\tau}_\mu=\operatorname{vec}(\tau_\mu)$ are the vectorized formats (through appropriate reshaping) of the far field, contrast,  and the unknown error respectively. The discrete inverse Born approximation then predicts the nodal values of $\mu$ by ignoring the error term $\vec{\tau}_\mu$ and using 
\begin{equation}
\label{reginv}
    \vec{\mu}_{\gamma}=(\mathcal{B}^{*}\mathcal{B}+\gamma I)^{-1}\mathcal{B}^*
\vec{u}_{\infty}\defeq \mathcal{B}_{\gamma}^{-1}\vec{u}_\infty,
\end{equation}
where  $\gamma>0$ is a regularization parameter. Experimentally, we observe that an optimal choice of $\gamma$ will generally be in the range $[0.1,1]$ for
the inverse problems in our study.

\subsection{Data generation}
\begin{figure}[t]
    \centering
    \includegraphics[width=0.7\linewidth]{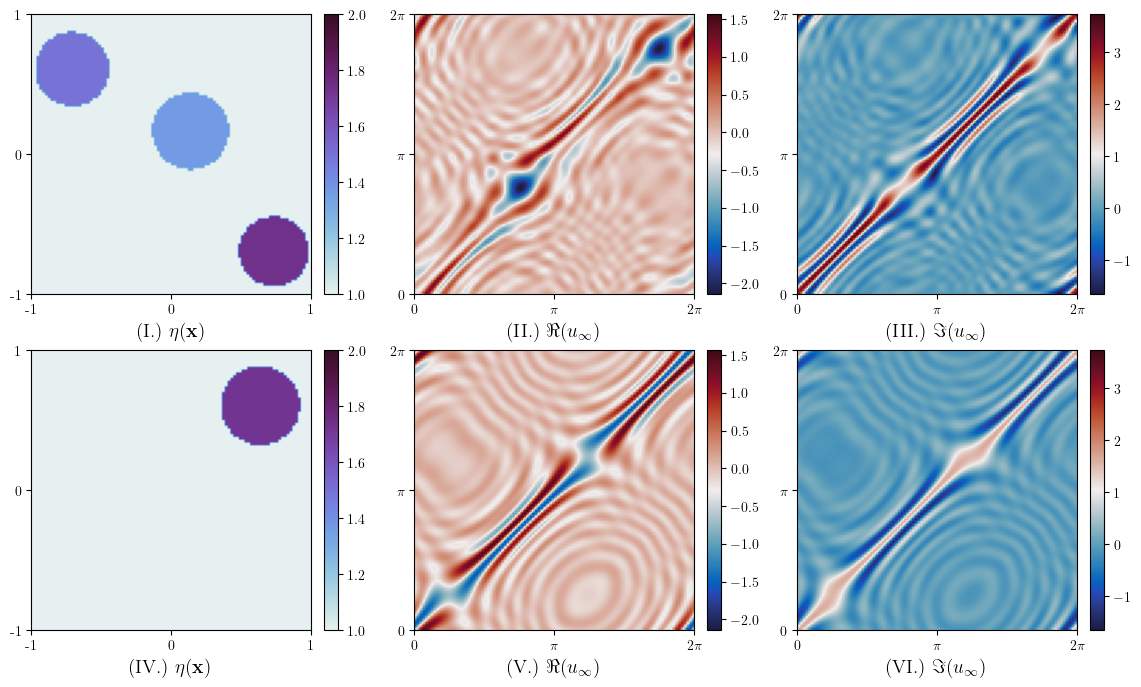}
    \caption{Two realizations of the random training data.  Each row displays a density plot of $\eta$ for a sample scatterer (left panel) alongside the real and imaginary components of its computed far field pattern.  The far field pattern is the data for the inverse scattering problem.}
    \label{fig:samples}
\end{figure}
As discussed at the end of Section \ref{sec:prob}, we utilize synthetic scattering data generated using a standard finite element approach to approximating (\ref{Helmholtz}). For our experiments, we fix the wave number at $k=16$ as in \cite{Fan:22} and set the number of sources and receivers to be $N_S=N_R=100$. Higher wave number problems may also be learned, provided $N_S$ and $N_R$ are increased. In the case of training data, we do not add extra measurement noise to the computed far field matrix $U_{\infty}$.
For a discussion of generalization to the case with noise, see Section \ref{noise}. For our training data, we also assume that no absorption occurs, i.e., $\Im(\eta(\textbf{x}))=0$. For a discussion of generalization to the case with weak absorption, see Section \ref{complex}. The spatial discretization of the domain $\Omega=[-1,1]^2$ uses $N_g=100$.

For training, we want to use simple scatterers. Here we create the scatterer field $\eta(\textbf{x})$ as the union of $N_{c}$ piecewise-constant circles overlaying a homogeneous background of air, where we uniformly sample $N_c\in \{1,2,3\}$. The $i$th circle is assigned a constant value $\eta_i$ sampled from the uniform distribution $\mathcal{U}(1.1, 1.8)$, with the background air set to $\eta=1$ (we shall see that $\eta\approx 1.8$ is outside the weak scattering regime). The radius and position of each circle are sampled uniformly from $\mathcal{U}(0.1,0.3)$ and $\mathcal{U}(-0.7, 0.7)$, respectively. In selecting the training data, for simplicity of mesh generation, we enforce that the circles cannot overlap and must be fully contained in the search domain $[-1,1]^2$. We generate 20,000 samples for training and validation (under an 80-20\% split) with an additional 4,000 samples for testing purposes. This choice is motivated by \cite{Fan:22} where 16,384 samples are used for training. Visualizations of typical training samples can be seen in Fig.~\ref{fig:samples}. 

The validity of the Born approximation is contingent upon weak scattering of the incident field. Moreover, the convergence of the Neumann series (and hence the accuracy of the Born approximation) requires $\eta$ to be sufficiently close to unity.   In the case of the synthetic dataset, the Born approximation fails when  $\eta_i\gg 1$ for any circle, thus  preventing an accurate inversion of (\ref{Tborn}). In particular, a direct inversion  results in severe underestimates of the true contrast regardless of the chosen regularization parameter $\gamma$.
When weak scattering breaks down, the Born approximation also gives rise to artifacts when there is multiple scattering.  
Examples of poor reconstructions using the Born approximation are shown in Section~\ref{sec:Num}.
We seek to remediate these failures in the strong scattering case through  corrective approaches using NNs.

\section{ Network architecture and training}\label{sec:NN}
Throughout the rest of the paper we shall used two norms defined on vectors with $n$ components so that if
$\vec{v}\in\mathbb{C}^n$ then the $\ell^p$ norm of $\vec{v}$ is
\[
\Vert \vec{v}\Vert_{p}=\left(\sum_{j=1}^n|v_j|^p\right)^{1/p},\;p=1,2.
\]
For arrays (in particular the array of pixel intensities for the image), we first reshape the array into a vector and then compute
the corresponding vector norm as above.  By abuse of notation we shall use $\Vert\cdot\Vert_p$ to indicate the norm of a vector or a matrix ($p=2$ is in fact the Frobenius norm of a matrix).
\subsection{Correction strategies}
We suggest two methods for correcting the Born approximation: a pre-correction and a post-correction. Let $\vec{\mu}_{\rm exact}$ denote the exact nodal values for the (vectorized) contrast. There is an error $\vec{\tau}_{\mu}$ associated with both the convergence of the Neumann series and its truncation (the Born approximation). We may  write 
\begin{equation}
    \label{CNNB}
    \mathcal{B}\vec{\mu}_{\rm exact}=\vec{u}_{\infty}+\vec{\tau}_{\mu}\defeq\vec{u}_{\infty}^{\tau_{\mu}}.
\end{equation}
where $\mathcal{B}$ is the stacked Born operator defined in (\ref{stackedBorn}), and $\vec{u}_\infty$ is the exact far field pattern. The factor $\vec{\tau}_{\mu}$ can be thought of as the discrepancy in the Born-obtained far field compared to the true  far field. A CNN is trained to predict $\vec{\tau}_{\mu}$, and this approach pre-corrects the far field data into a form suitable for the Born approximation to predict an accurate contrast. In particular 
$$\mathcal{B}^{-1}_{\gamma}(\vec{u}_{\infty}^{\tau_{\mu}})=\mathcal{B}_{\gamma}^{-1}(\mathcal{B}\vec{\mu}_{\rm exact})\approx \vec{\mu}_{\rm exact}$$
provided that $\gamma$ is chosen correctly and $\vec{u}_{\infty}^{\tau_{\mu}}$ is learned appropriately. This is motivated by the approach in \cite{koponen}, where a NN is employed to compensate for modeling errors introduced by approximate forward models.

The second strategy to correct the Born approximation is to obtain an accurate representation of $\varepsilon_{\mu}$, the  discrepancy between the naive Born reconstruction and the true contrast. The factor $\varepsilon_{\mu}$ encompasses the far field behavior excluded by the Born approximation as well as the numerical error from the inversion scheme. We write
\begin{equation}
\varepsilon_{\mu}= \vec{\mu}_{\rm exact}-\mathcal{B}_{\gamma}^{-1}\vec{u}_{\infty},
\end{equation}
where, as in (\ref{reginv}), ${\mathcal B}_{\gamma}^{-1}$ denotes the regularized inverse of the Born operator. Once $\varepsilon_{\mu}$ is learned, an improved contrast estimate can be obtained through the simple formula
$$\vec{\mu}_{\rm exact}\approx \mathcal{B}^{-1}_{\gamma}\vec{u}_{\infty}+\mathcal{\varepsilon}_{\mu}.$$
This is the approach used in, for example, \cite{Taskinen:22,Lipponen:22} in a different context for correcting satellite data-based retrieval algorithms. 

\subsection{Training}\label{NIO}
\begin{figure}[t]
    \centering
   \begin{tikzpicture}[
    node distance=1cm and 0.5cm,
    box/.style={draw, rounded corners, minimum width=1.5cm, minimum height=1cm, align=center},
    arrow/.style={-Latex, thick},
    dashedbox/.style={draw=black, dashed, rounded corners, inner sep=0.4cm}
    ]

\node[box, fill=cyan!20] (ReU3) {$\Re(\vec{u}_{\infty})$};
\node[box, fill=cyan!20, below=0.1cm of ReU3] (ImU3) {$\Im(\vec{u}_{\infty})$};
\node[box, fill=red!20, right=1.5cm of ReU3.east|-ImU3.east, yshift=0.5cm] (NN3) {CNN/NIO/BCR};
\node[box, fill=green!20, right=of NN3] (M3) {$\eta(\textbf{x})$};

\draw[arrow] (ReU3) -- (NN3);
\draw[arrow] (ImU3) -- (NN3);
\draw[arrow] (NN3) -- (M3);

\node[dashedbox, fit=(ReU3) (ImU3) (NN3) (M3)] (box1) {};
\node[left=0.5cm of box1.west] {(a)};

\node[box, fill=cyan!20, below=1cm of ImU3] (ReU2) {$\Re(\vec{u}_{\infty})$};
\node[right=-0.24cm of ReU2](RightBorderReU2){};
\node[box, fill=cyan!20, below=1cm of ReU2] (ImU2) {$\Im(\vec{u}_{\infty})$};
\node[right=-0.24cm of ImU2](RightBorderImU2){};
\node[box, fill=blue!20, right=1.5cm of ReU2.east|-ImU2.east, yshift=2cm] (B2) {$\mathcal{B}_{\gamma}^{-1}$};
\node[box, fill=pink!20, right=of B2] (TildeM) {$\tilde{\eta}(\textbf{x})$};
\node[box, fill=red!20, below right=of ReU2, xshift=1cm] (NN2) {BCNN};
\node[box, fill=yellow!20, right=of NN2] (Epsilon) {$\varepsilon$};
\node[box, fill=green!20, right=of Epsilon] (M2) {$\eta(\textbf{x})=\tilde{\eta}(\textbf{x})+\varepsilon$};

\draw[arrow] (RightBorderReU2) -- (B2);
\draw[arrow] (RightBorderImU2) -- (B2);
\draw[arrow] (RightBorderReU2) -- (NN2);
\draw[arrow] (RightBorderImU2) -- (NN2);
\draw[arrow] (B2) -- (TildeM);
\draw[arrow] (TildeM) -- (M2);
\draw[arrow] (NN2) -- (Epsilon);
\draw[arrow] (Epsilon) -- (M2);

\node[dashedbox, fit=(ReU2) (ImU2) (B2) (TildeM) (NN2) (Epsilon) (M2)] (box2) {};
\node[left=0.5cm of box2.west] {(b)};

\node[box, fill=cyan!20, below=1cm of ImU2] (ReU1) {$\Re(\vec{u}_{\infty})$};
\node[box, fill=cyan!20, below=0.1cm of ReU1] (ImU1) {$\Im(\vec{u}_{\infty})$};
\node[box, fill=red!20, right=1.5cm of ReU1.east|-ImU1.east, yshift=0.5cm] (NN1) {CNNB};
\node[box, fill=purple!20, right=of NN1] (VF1) {$\vec{u}_{\infty}^{\tau}$};
\node[box, fill=blue!20, right=of VF1] (B1) {$\mathcal{B}_{\gamma}^{-1}$};
\node[box, fill=green!20, right=of B1] (M1) {$\eta(\textbf{x})$};

\draw[arrow] (ReU1) -- (NN1);
\draw[arrow] (ImU1) -- (NN1);
\draw[arrow] (NN1) -- (VF1);
\draw[arrow] (VF1) -- (B1);
\draw[arrow] (B1) -- (M1);

\node[dashedbox, fit=(ReU1) (ImU1) (NN1) (VF1) (B1) (M1)] (box3) {};
\node[left=0.5cm of box3.west] {(c)};

\end{tikzpicture}
    \caption{The neural network architecture for the three models compared in this paper. (a) The CNN, NIO, and BCR models directly map the far field to the squared refractive index. Note that the BCR model only uses $\Re(\vec{u}_{\infty})$ as input. (b) The regularized Born inverse is computed to obtain a rough estimate of the squared refractive index $\eta(\textbf{x})$, which is then corrected by the BCNN predicted $\varepsilon$. (c) The CNNB model pre-corrects the far field and then applies a regularized Born inverse to obtain the squared refractive index. Note that the necessary reshaping of the far field and contrast are omitted in this diagram, as well as the map from the contrast to the refractive index: $\mu\mapsto \mu+1=\eta$.}
    \label{fig:architecture}
\end{figure}
For the two correction strategies, we can associate the following CNN-Born (CNNB) and Born-CNN (BCNN) models, respectively (see Fig.~\ref{fig:architecture}). The former performs the Born approximation on the corrected input far field while the latter performs the Born approximation on the labels of the training set. For a training set of size $N_{\rm train}$ we have the corresponding training regimes:
\begin{itemize}
    \item \textbf{CNNB}: We consider data pairs $\{(\vec{u}_{\infty}^{(i)},\vec{u}^{\tau_{\rev{\mu}} (i)}_{\infty})\}_{i=1}^{N_{\rm train}}$ and seek to learn $\vec{\tau}_\mu$ by minimizing the loss function
\begin{equation}
\label{lossCNNB}
    \mathcal{L}_{\rm CNNB}(\beta)=\frac{1}{N_{\rm train}}\sum_{i=1}^{N_{\rm train}}\|\vec{u}^{\tau_{\mu}(i)}_{\infty}-\operatorname{CNNB}(\vec{u}_{\infty}^{(i)};\beta)\|_2^2,
\end{equation}
where $\rm CNNB(\:\cdot\:;\beta)$ is a CNN with weights and biases collected in the vector  $\beta$.
\item \textbf{BCNN}: We consider data pairs $\{(\vec{u}_{\infty}^{(i)},\varepsilon_{\mu}^{(i))}\}_{i=1}^{N_{\rm train}}$ and seek to learn $\vec{\varepsilon}_\mu$ by minimizing the loss function 
\begin{equation}
\label{lossBCNN}
    \mathcal{L}_{\rm BCNN}(\beta')=\frac{1}{N_{\rm train}}\sum_{i=1}^{N_{\rm train}}\|\varepsilon_{\mu}^{(i)}-\operatorname{BCNN}(\vec{u}_{\infty}^{(i)};\beta')\|_2^2,
\end{equation}
where $\rm BCNN(\:\cdot\:;\beta')$ is another CNN with weights and biases collected in the vector $\beta'$.
\end{itemize}

To provide a baseline for comparison, we train a simple black-box CNN to directly map the far field data to the solution of the inverse scattering problem motivated by
\cite{Fan:22}. In other words, we consider data-pairs $\{(\vec{u}_{\infty}^{(i)},\vec{\mu}^{(i)})\}_{i=1}^{N_{\rm train}}$ and seek to directly learn $\vec{\mu}$ by minimizing the loss function
\begin{equation}
\label{lossCNN}
    \mathcal{L}_{\rm CNN}(\beta'')=\frac{1}{N_{\rm train}}\sum_{i=1}^{N_{\rm train}}
\|\vec{\mu}^{(i)}-{\rm CNN}(u_{\infty}^{(i)};\beta'')\|_2^2,
\end{equation}
where $\rm CNN(\:\cdot\:,\beta'')$ is a third CNN with weights and biases $\beta''$. 
To aid comparison we use the same general CNN structure in all three cases, but optimize the network hyper-parameters to each case. 
Furthermore, we compute the basic regularized Born inverse (\ref{reginv}) with regularization parameter $\gamma=1$ and $\gamma=0.1$ for a direct comparison of CNN results to the Born approximation.

\begin{table}[t]
\caption{Hyper-parameters selected by tuning the CNN, BCNN, and CNNB models.}
\label{table:opt}
\centering
\begin{tabular}{c|| c c c}
\hline
 & CNN & BCNN & CNNB \\ [0.5ex]
\hline\hline
Conv2D Layers & 4 & 4 & 4 \\ 
Conv2D Channels & [296, 211, 152, 61] & [335, 33, 195, 65] & [125, 358, 426, 221] \\ 
FC Layers & 3 & 1 & 1 \\ 
FC Units  & [537, 465, 419] & [971] & [576] \\ 
Activation & GELU & GELU & GELU \\ [1ex]
\hline
\end{tabular}
\end{table}
Besides the straight-foward CNN, we also compare to results computed using the implementation of NIO~\cite{NIO} available from Github\footnote{\href{https://github.com/mroberto166/nio}{https://github.com/mroberto166/nio}}.  Although there is a Helmholtz-based example in \cite{NIO}, it is not for scattering on an infinite domain (rather it is for Dirichlet and Neumann data pairs on a bounded domain). Nonetheless, the software can  be easily modified to handle far field data from incident plane waves as in our case.  We trained  NIO
using our data as discussed previously for CNN with an $\ell_2$ loss function.  Using scripts from the GitHub repository, we optimized NIO for our problem.

We also compare to the BCR-Net  of \cite{FAN20191,Fan:22}, trained using our data with $\ell^2$ loss and utilizing the training scripts and network code generously provided to us by its authors.  Training on our data is necessary to make a comparison with our networks meaningful.
In the CNN and NIO approaches, the model takes as input a representation of the far field matrix,  where the real and imaginary parts are separated and stored in two distinct channels, which enables the model to process both components simultaneously. Meanwhile, the BCR-Net model takes solely the real part of the far field matrix as input. The model output is then reshaped into a two-dimensional map corresponding to the discrepancy term. This end-to-end mapping allows the model to leverage spatial correlations in the data that are challenging to capture through purely analytic approaches. 
The high-level structure of the CNN based methods is visualized in Fig.\ref{fig:architecture}. For the structure of BCR-Net see~\cite{Fan:22} and for NIO see~\cite{NIO}.

The architecture of each CNN was determined through a randomized search process, where 150 network configurations were generated, and the best-performing model was selected based on the lowest validation mean squared error (MSE). The randomization process involved varying hyper-parameters, including the number of convolution layers (ranging from 1 to 4) and the number of fully connected (FC) layers (ranging from 1 to 3). The number of channels per each convolution layer was randomly chosen between 16 and 512, while the number of units per FC layer was selected between 64 and 1024. Max pooling, with a fixed kernel size of 2, was applied after each convolutional layer.
For activation functions, we randomly selected from a set including Rectified Linear Unit (ReLU), LeakyReLU with a negative slope of 0.1,  Gaussian Error Linear Unit (GELU), Sine, and Sigmoid. The last layer always used a linear activation function. All CNN based models were trained using the AdamW optimizer (discussed in \cite{loshchilov2019decoupledweightdecayregularization}) with an initial 
learning rate of 0.0005, adjusted via a learning rate scheduler that applied a 5\% decay every 100 epochs. The models were trained with an early stopping patience of 50 epochs. 
Table \ref{table:opt} provides a summary of the tuned hyper-parameters.

\section{Numerical experiments}\label{sec:Num}
\subsection{In-distribution performance}
Following the methods discussed in Section \ref{sec:prob}, we generated 4,000 in-distribution scatterers to evaluate model performance. 
Table \ref{table:l2} summarizes the average $\ell^2$ and $\ell^1$ errors of the test scatterer reconstructions.
Table~\ref{table:l2} shows that while NIO and BCR-Net slightly outperform the CNN based models in both norms, all NN based models perform similarly well in $\ell^2$ (training used this norm for the loss function).
The results also show that all NN based models significantly outperform the Born approximation in terms of accuracy, reducing the relative $\ell^2$ error by almost 50\%. This is likely because the Born approximation performs poorly for large contrast scatterers which are included in the  testing data. 

\begin{figure}[H]
    \centering
    \includegraphics[width=\linewidth]{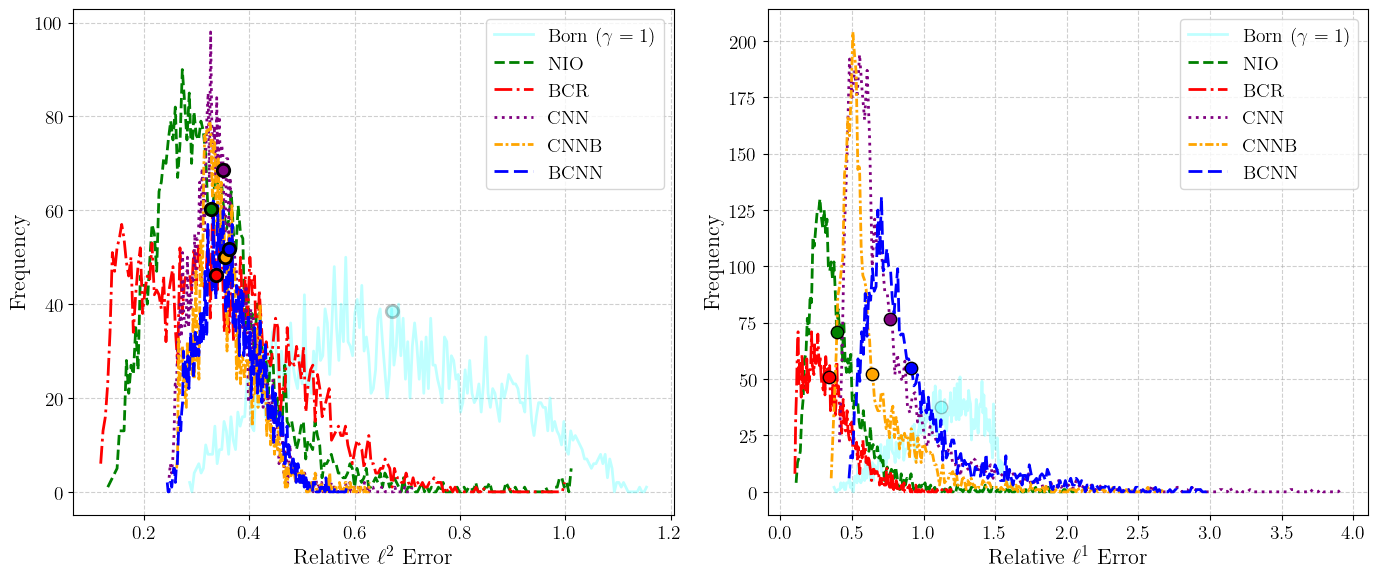}
    \caption{The distribution of relative $\ell^2$ (left) and $\ell^1$ (right) norm errors for the test dataset. The average of each distribution is denoted by the corresponding colored circle. The NN based methods all improve over the 
    simple Born approximation. }
    \label{fig:error}
\end{figure}

\begin{figure}[H]
    \centering
    \includegraphics[width=0.7\linewidth]{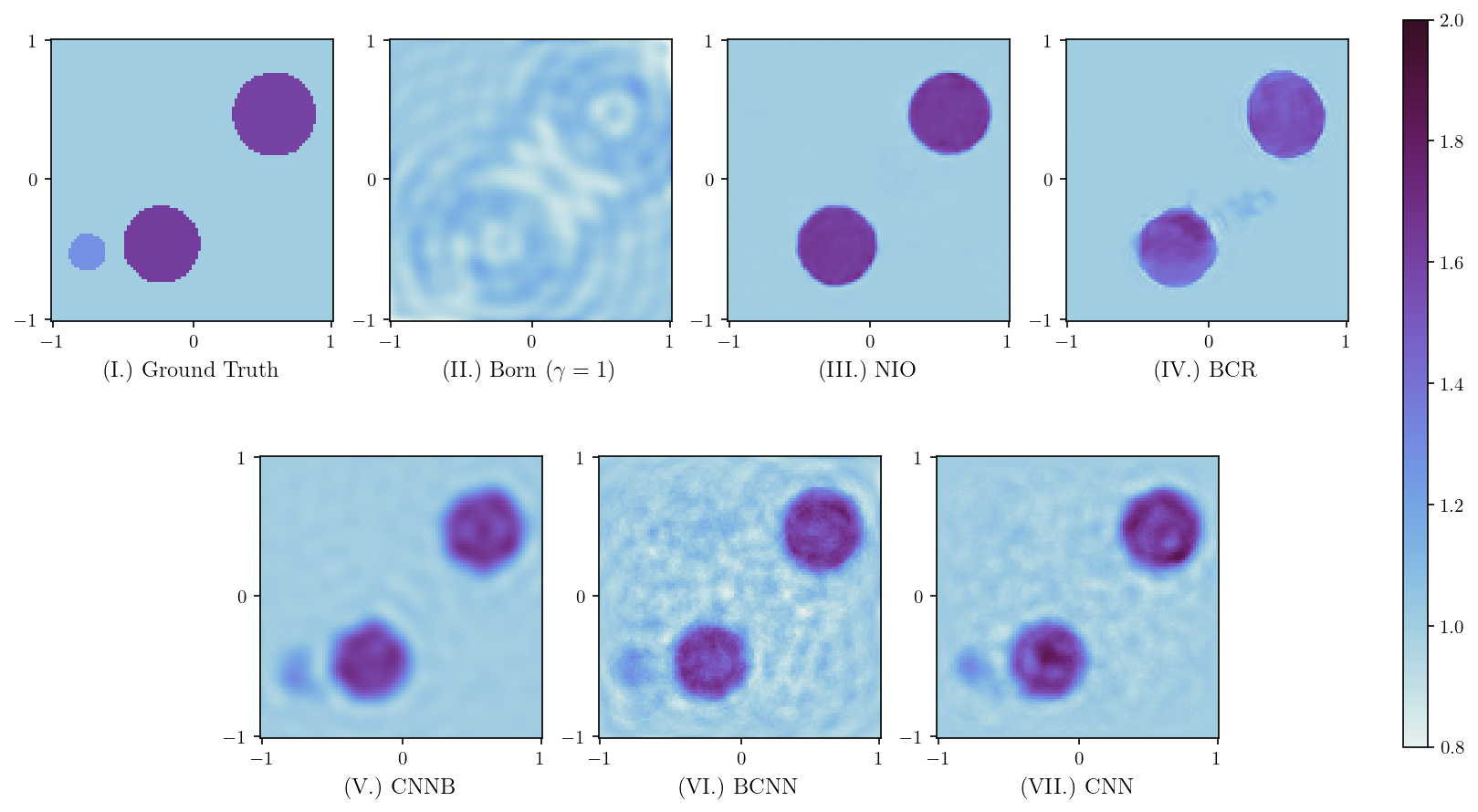}
    \caption{A sample reconstruction of three in distribution circular scatterers using each model. All NN based methods are much more accurate than the Born approximation since, in this example, the scattering is strong. While NIO (III) and BCR (IV) produce sharper edges of the scatterer, in this example they completely fail to detect the smaller weak scatterer.}
    \label{fig:test_res}
\end{figure}

To investigate further how the NN based models compare to the inverse Born approximation by itself, Fig.~\ref{fig:error} visualizes the error distribution for the test samples. The results show that all NN based models significantly outperform the basic inverse Born approximation in terms of accuracy, reducing the relative $\ell^2$ error by almost 50\%. This is likely because the Born approximation performs poorly for large contrast scatterers which are included in the  testing data. 
Fig~\ref{fig:test_res}  shows the actual profile  of one test scatterer together with the reconstructions.
Qualitatively, we  observe that NIO and BCR-Net sometimes  have a tendency to forgot smaller, weaker scatterers.
Going forward, we will omit the error and reconstruction results related to the inverse Born approximation for $\gamma=0.1$ due to its consistently sub-par performance in relation to when $\gamma=1$.

\begin{table}[t]
\caption{Average relative $\ell^2$ and $\ell^1$ error (\%) in the reconstruction of the contrast $\mu$ on the test dataset of $4,000$ scatterers. The best reconstruction for each metric is in boldface. \label{table:l2}}
\centering 
\begin{tabular}{|c|| c c|}
\hline
Method & $\ell^2$ (\%) & $\ell^1$ (\%)\\ [0.5ex]
\hline\hline
Born  ($\gamma=1$) & 67.0378 & 111.8617 \\
Born ($\gamma=0.1)$ & 74.5836 & 152.5471 \\
NIO & \textbf{32.6958} & 39.2627 \\
BCR-Net & 33.6819 & \textbf{33.6776}\\
CNN  & 35.0081 & 76.5449\\
CNNB & 35.3143 & 63.9699 \\
BCNN & 36.1794 & 91.2601 \\
[0.5ex]
\hline\end{tabular}
\end{table}

\subsection{Noise robustness}
\label{noise}

\begin{figure}[H]
    \centering
    \includegraphics[width=\linewidth]{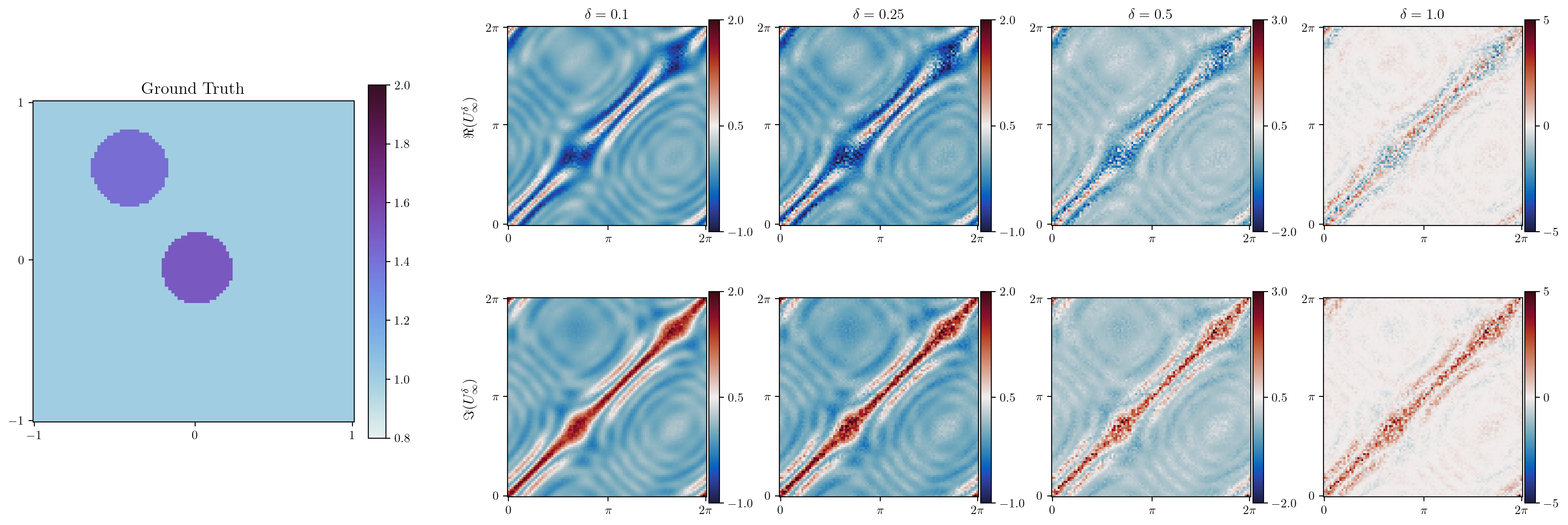}
    \caption{Visualization of a scatterer and its noisy far-field patterns for multiple noise levels in correspondence with (\ref{farfieldmat}). The left column shows the refractive index of the scatterer. The top and bottom rows display the real and imaginary parts of the far-field data 
    for $\delta=0.1,0.25,0.5,1.0$, respectively. The reconstruction of this scatterer is shown in Fig.~\ref{fig:noise_predict}.}
    \label{fig:noise_problem}
\end{figure}

Noise, arising from practical measurement apparatus errors or environmental factors, often corrupts scattering data, posing significant challenges to accurately reconstructing scatterers. We trained each model using noise free data (apart from numerical error), so it is important to determine how the various models cope with out-of-distribution data corresponding to noisy far field data. In this context, following \cite{Fan:22},  we consider a far field matrix $U^\delta_{\infty}$ affected by noise, with entries given by
\begin{equation}
\label{farfieldmat}
    (U^\delta_{\infty})_{i,j}\defeq (1+\delta \mathcal{E}_{i,j})u_{\infty}(\hat{\textbf{x}}_i,\textbf{d}_j),
\end{equation}
where $\delta$ controls the strength of the noise and $\mathcal{E}_{i,j}$ is the (complex-valued) noise for incoming wave direction $\varphi_i$ and measurement direction $\theta_j$. To model the noise, we assume that $\mathcal{E}_{i,j}$ is sampled from a uni-variate complex standard normal distribution ${\cal CN}(0,1)$ \cite[Def 2.1]{Andersen1995}.

\begin{table}[H]
\caption{Average relative $\ell^2$ and $\ell^1$ errors (\%) in the reconstruction of the contrast $\mu$ on the test dataset of 4000 scatterers with noise added to the  far field data having parameter $\delta=0.1,0.25,0.5,1.0$ defined in (\ref{farfieldmat}). The best reconstruction for each $\delta$ is in boldface. }
\centering
\renewcommand{\arraystretch}{1.2}
\begin{tabular}{|c|c c||c|c c|}
\hline
\multicolumn{3}{|c||}{\textbf{$\delta=0.1$}} & 
\multicolumn{3}{c|}{\textbf{$\delta=0.25$}} \\
\hline
Method & $\ell^2$ (\%) & $\ell^1$ (\%) & Method & $\ell^2$ (\%) & $\ell^1$ (\%) \\
\hline
Born ($\gamma=1$)   & 67.0378 & 111.8617 & Born ($\gamma=1$)   & 67.0378 & 111.8617 \\
NIO                 & 36.8788 & 45.1435 & NIO                 & 63.1058 &  160.1682\\
BCR-Net &\textbf{33.8625}& \textbf{33.8459}&BCR-Net&\textbf{34.8135} &\textbf{34.8141}\\
CNN                 & 36.1031 & 78.3238 & CNN                 & 41.9604 &  90.9867\\
CNNB                & 35.8799 & 65.9740 & CNNB                & 39.2874 &  79.4563\\
BCNN                & 38.4137 & 102.7461 & BCNN                & 48.4455 & 143.9246 \\
\hline
\multicolumn{3}{|c||}{\textbf{$\delta=0.5$}} & 
\multicolumn{3}{c|}{\textbf{$\delta=1.0$}} \\
\hline
Method & $\ell^2$ (\%) & $\ell^1$ (\%) & Method & $\ell^2$ (\%) & $\ell^1$ (\%) \\
\hline
Born ($\gamma=1$)   & 67.0378 & 111.8617 & Born ($\gamma=1$)   & 67.0378 & 111.8617 \\
NIO                 & 103.4922 & 358.3474 & NIO                 & 148.9648 & 535.9080 \\
BCR-Net &\textbf{37.8709} &\textbf{38.2524} &BCR-Net &\textbf{47.3537} &\textbf{51.6489}\\
CNN                 & 62.6698 & 135.3118 & CNN                 & 134.2359 & 354.2700 \\
CNNB                & 49.2495 & 116.5801 & CNNB                & 73.5768 & 183.3327 \\
BCNN                & 89.3354 & 281.7972 & BCNN                & 380.3745 & 1173.8297 \\
\hline
\end{tabular}
\label{noise_tab}
\end{table}
We analyze the zero-shot performance of the models on the testing dataset described previously under varying levels of added random noise using $\delta=0.1$, $0.25$, $0.5$, $1$. We also compare their reconstructions to the Born model ($\gamma=1$), which is particularly robust to noise due to the regularized inversion, yet under-approximates the contrast for strong scatterers. 
Table~\ref{noise_tab} provides statistics on the behavior of the models.

\begin{figure}[H]
    \centering
    \includegraphics[width=0.9\linewidth]{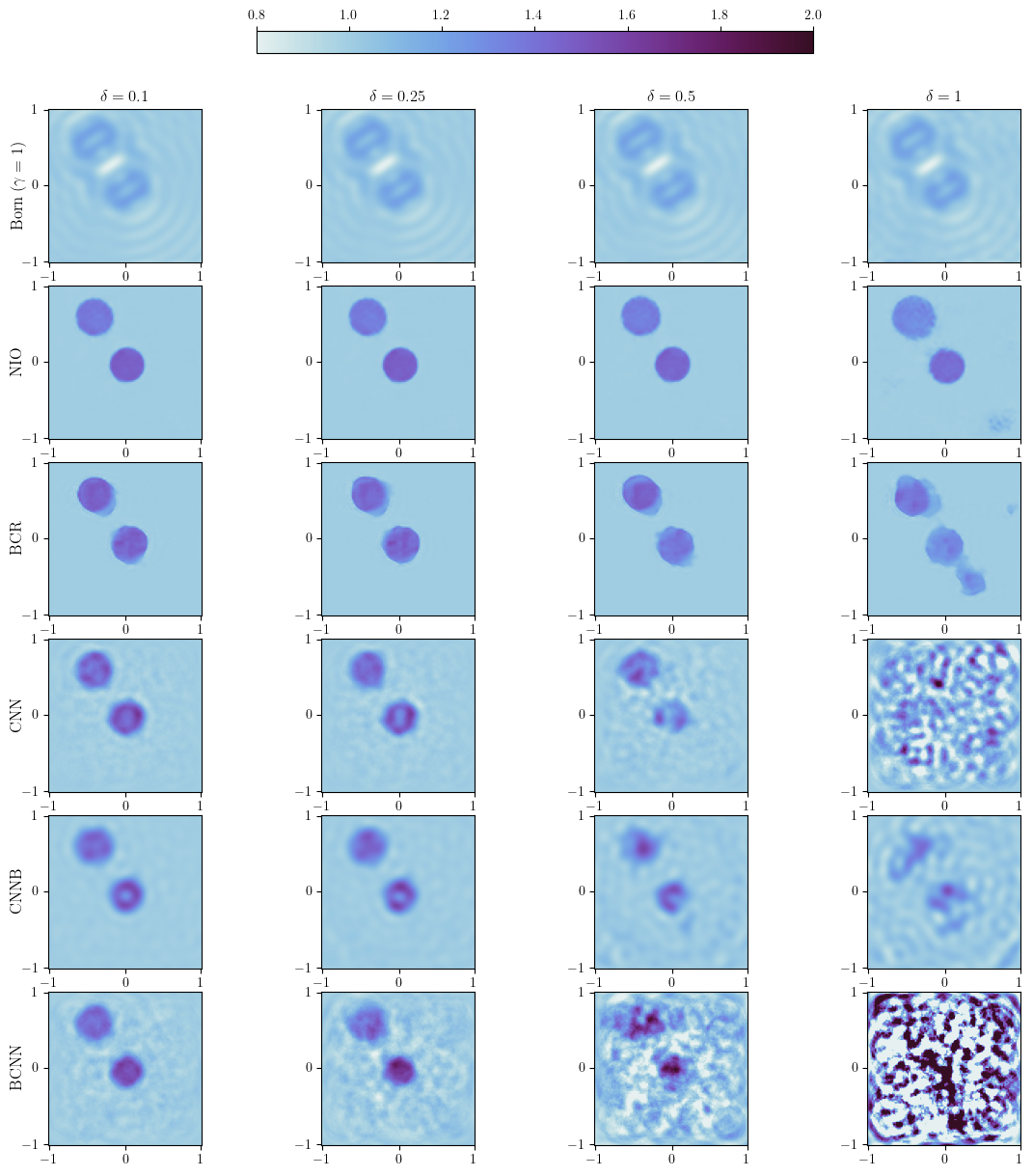}
    \caption{The refractive index estimated by each model corresponding to the scatterer in Fig.~\ref{fig:noise_problem} with noisy data. In this example, the high noise regime, NIO and BCR tend to hallucinate additional scatterers while BCNN and CNN are marred with artifacts. Meanwhile, CNNB demonstrates relative robustness even at $\delta=1$ without becoming highly corrupted by background artifacts. }
    \label{fig:noise_predict}
\end{figure}

At all levels of noise, Table~\ref{noise_tab} shows that
BCR-Net is the most stable to noise. All NN based schemes 
perform roughly equivalently for lower noise levels at $\delta=0.1,0.25$ with CNNB being the most stable of the Born based techniques. 
The standard Born approximation with $\gamma=1$ is strongly regularized to noise, but cannot handle large contrasts resulting in poor average performance.
We conjecture that the Born inversion incorporated after the CNN in CNNB  smooths high frequency error effectively controlling noise amplification and making it the more stable of the Born based techniques.  In contrast, BCNN first provides a smooth reconstruction, and then applies the CNN using that reconstruction and the far field data to derive a corrected model. This could result in increased sensitivity to noise since CNN is used after regularization.
Fig.~\ref{fig:noise_predict} shows the results for 
one reconstruction at various noise levels.  In this case the ground truth scatterer is shown in Fig.~\ref{fig:noise_problem}. At high noise levels, this example shows that both NIO and BCR-Net may hallucinate as might be expected~\cite{siam-rev}. BCNN and CNN do not perform as well as CNNB at the highest noise level.

\newpage
\subsection{Absorption}
\label{complex}

\begin{figure}[H]
    \centering
    \includegraphics[width=0.7\linewidth]{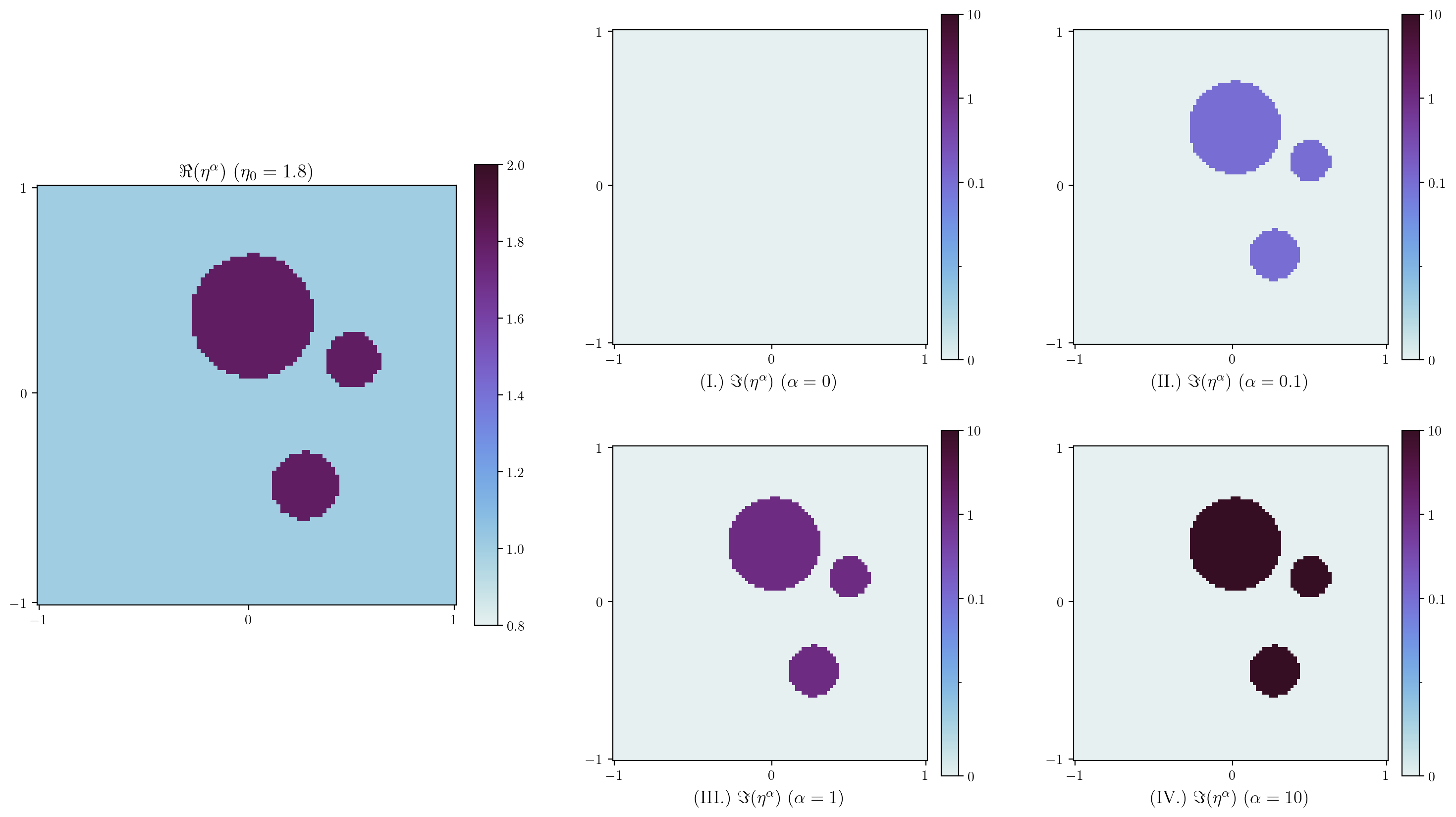}
    \caption{Visualization of an absorbing scatterer for multiple absorption levels in correspondence with \ref{absorbeta}. The left column shows the real part of the scatterer with $\eta_0=1.8$. The right column displays the imaginary part of the scatterer for $\alpha= 0, 0.1, 1.0, 10.0$, respectively. The reconstruction of the scatterer is shown in Fig.~\ref{fig:complex_corrupt}}.
    \label{fig:absorption_problem}
\end{figure}

Our models are trained on data from absorption free media
($\Im(\eta)=0$), as would be appropriate for applications
involving non-dissipative media. In many practical applications, the media exhibits varying degrees of absorption~\cite{Jackson}, characterized by $\Im(\eta) > 0$. Of course, the inverse Born approximation handles absorbing  media seamlessly (provided the scattering is weak), but the NN based schemes depend on the training data.

Since we did not train the NN based models using far field data from absorbing scatterers, it is important to determine if the NN based models can predict the real part of $\eta$ even when the scattering media are slightly absorbing.    

\begin{table}[h]
\caption{Average relative $\ell^2$ and $\ell^1$ errors (\%) in the reconstruction of the contrast $\mu$ on the test dataset of 100 absorbing scatterers with $\alpha=0,0.1,1,10$ following (\ref{absorbeta}). The best reconstruction for each $\alpha$ is in boldface.}
\centering
\renewcommand{\arraystretch}{1.2}
\begin{tabular}{|c|c c||c|c c|}
\hline
\multicolumn{3}{|c||}{\textbf{$\alpha=0$}} & 
\multicolumn{3}{c|}{\textbf{$\alpha=0.1$}} \\
\hline
Method & $\ell^2$ (\%) & $\ell^1$ (\%) & Method & $\ell^2$ (\%) & $\ell^1$ (\%) \\
\hline
Born ($\gamma=1$)   & 88.9285 & 144.9036 & Born ($\gamma=1$)   & 89.4273 & 133.2401 \\
NIO                 & \textbf{32.5428} & \textbf{36.7441}  & NIO                 & \textbf{33.3716} & \textbf{38.7768} \\
BCR-Net             & 41.2125     &45.9028          &  BCR-Net &42.7416 &   49.2050      \\  
CNN                 & 37.7394 & 78.1480  & CNN                 & 39.9897 & 77.3409 \\
CNNB                & 37.4499 & 65.0130  & CNNB                & 38.2855 & 68.8976 \\
BCNN                & 41.6169 & 107.8156 & BCNN                & 43.1420  & 103.9886 \\
\hline
\multicolumn{3}{|c||}{\textbf{$\alpha=1$}} & 
\multicolumn{3}{c|}{\textbf{$\alpha=10$}} \\
\hline
Method & $\ell^2$ (\%) & $\ell^1$ (\%) & Method & $\ell^2$ (\%) & $\ell^1$ (\%) \\
\hline
Born ($\gamma=1$)   & 96.2615 & 119.4021 & Born ($\gamma=1$)   & 95.1858 & 147.9940 \\
NIO                 & 67.1786 & \textbf{83.2703}  & NIO                 & 94.7267 & 144.64945 \\
BCR-Net             &\textbf{61.9917}       &  85.0740       &  BCR-Net & \textbf{72.26803}        &    \textbf{117.4409}      \\  
CNN                 & 71.2832 & 109.7257 & CNN                 & 78.2193 & 137.1136 \\
CNNB                & 70.4011 & 107.2159 & CNNB                & 77.0885 & 136.1209 \\
BCNN                & 70.9577 & 139.9843 & BCNN                & 78.3788 & 188.4496 \\
\hline
\end{tabular}
\label{absorb_tab}
\end{table}

To investigate this scenario, we consider the case in which $\eta$ inside a circular scatterer is perturbed by an unknown, small absorption term. Specifically, we use the following complex-valued $\eta$ to generate far field test data:
\begin{equation}
\label{absorbeta}
    \eta^{\alpha} \defeq \eta_0+\iu\alpha,
\end{equation}
where $\alpha > 0$ controls the strength of the perturbation and $\eta_0$ is the real part of the refractive index of the scattering medium. For this experiment, we fix $\eta_0=1.8$ in the scatterer. Since we assume that the background medium is air that is not absorbing, we do not add any absorption to the background, which remains at $\eta=1$. We specifically consider the cases where $\alpha=0,\:1/10,\:1,\:10$.  

The results are shown in Table~\ref{absorb_tab}. We observe that, judged using the $\ell^2$ norm, the NIO architecture performs best under weak absorption for $\alpha=0.1$ in the scatterer. CNNB also tolerates small absorption and is the best of the CNN based models. However, all models
deteriorate  for higher magnitudes of absorption ($\alpha=1,10$) and BCR-Net  performs optimally in the $\ell^2$ in this regime. 

\begin{figure}[H]
    \centering
    \includegraphics[width=0.9\linewidth]{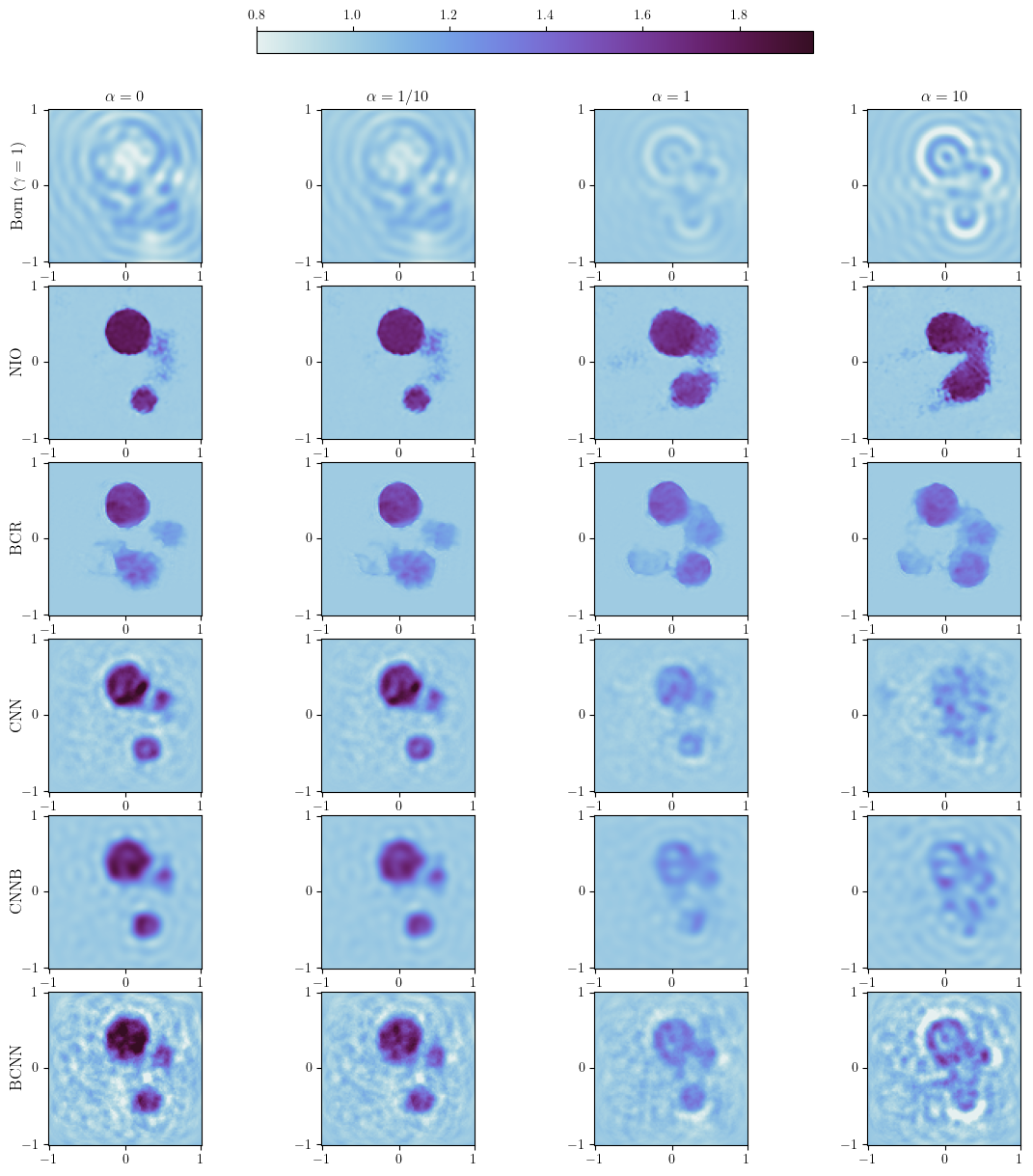}
    \caption{
    Here we show the refractive index estimated by each model with data from an absorbing scatterer corresponding to the ground truth  scatterer ($\eta_0=1.8$) in Fig.~\ref{fig:absorption_problem}. For larger absorption values, the CNN based models tend to under-approximate with $\eta_0\approx 1.3$ and produce blurry artifacts surrounding the scatterers. Meanwhile, NIO over-estimates with $\eta_0\approx 2.1$ and fails to distinguish the smaller scatterer, while BCR-Net predicts a spurious fourth scatterer. }
    \label{fig:complex_corrupt}
\end{figure}

Fig.~\ref{fig:complex_corrupt} shows the results for one reconstruction at various absorption levels of the ground truth scatterer in Fig.~\ref{fig:absorption_problem}. In this case, NIO does not successfully reconstruct the smallest scatterer. BCR-Net places this scatterer further from the large scatterer and hallucinates by predicting a fourth scatterer for $\alpha=1,10$. 

\newpage
\subsection{Increased scatterer complexity}\label{sec:isc}
We next evaluate the generalization capabilities of the models when applied to more complex out-of-distribution scatterers. For each scatterer type, we generate 100 instances with random type-specific variations. The performance of the various reconstruction algorithms are summarized in Table~\ref{tab:errors} and sample instances are displayed in Fig.~\ref{fig:u_shape}--\ref{fig:squares}.

\begin{table}[h]
\caption{Average relative $\ell^2$ and $\ell^1$ errors (\%) in the
reconstruction of the contrast $\mu$ on 100 randomized samples of the complex
scatterers discussed in Section~\ref{sec:iscat}. The best reconstruction for
each problem is in boldface.}
\centering
\renewcommand{\arraystretch}{1.2}
\begin{tabular}{|c|c c||c|c c|}
\hline
\multicolumn{3}{|c||}{\textbf{U}} & 
\multicolumn{3}{c|}{\textbf{Ring}} \\
\hline
Method & $\ell^2$ (\%) & $\ell^1$ (\%) & Method & $\ell^2$ (\%) & $\ell^1$ (\%) \\
\hline
Born ($\gamma=1$)   & 69.6645 & 119.4050 & Born ($\gamma=1$)   & 109.3183 & 116.4062 \\
CNN                 &  55.8611 & 110.5492 & CNN                 & 64.7712 & 65.2693 \\
NIO                 & 77.8722 & 109.2033 & NIO                 & 94.3889 & 98.3254 \\
CNNB                & \textbf{53.5359} & \textbf{98.3680} & CNNB                & \textbf{60.5713} & \textbf{63.8227} \\
BCNN                & 53.9426 & 130.5083 & BCNN                & 74.5521 & 82.6612 \\
BCR-Net             &    80.3159    &  105.6085  & BCR-Net &  73.9569        & 74.0454  \\
\hline
\multicolumn{3}{|c||}{\textbf{Shepp-Logan}} & 
\multicolumn{3}{c|}{\textbf{Rectangles}} \\
\hline
Method & $\ell^2$ (\%) & $\ell^1$ (\%) & Method & $\ell^2$ (\%) & $\ell^1$ (\%) \\
\hline
Born ($\gamma=1$)   & 97.9434 & 138.3105 & Born ($\gamma=1$)   & 89.5299 & 118.3325 \\
CNN                 & 59.2884 & \textbf{84.0929} & CNN                 & 42.6187 & 58.0487 \\
NIO                 & 75.05212 & 91.7023 & NIO                 & 57.4389 & 59.3596 \\
CNNB                & \textbf{58.5902} & 87.1067 & CNNB                & \textbf{42.0379} & \textbf{55.1703} \\
BCNN                & 59.0327 & 102.3163 & BCNN                & 45.9519 & 75.3718 \\
BCR-Net             &  78.3657  & 87.3097        & BCR-Net & 73.8079        &     75.08584     \\
\hline
\end{tabular}
\label{tab:errors}
\end{table}

\subsubsection{Resonance Structure}
\begin{figure}[t]
    \centering
    \includegraphics[width=0.7\linewidth]{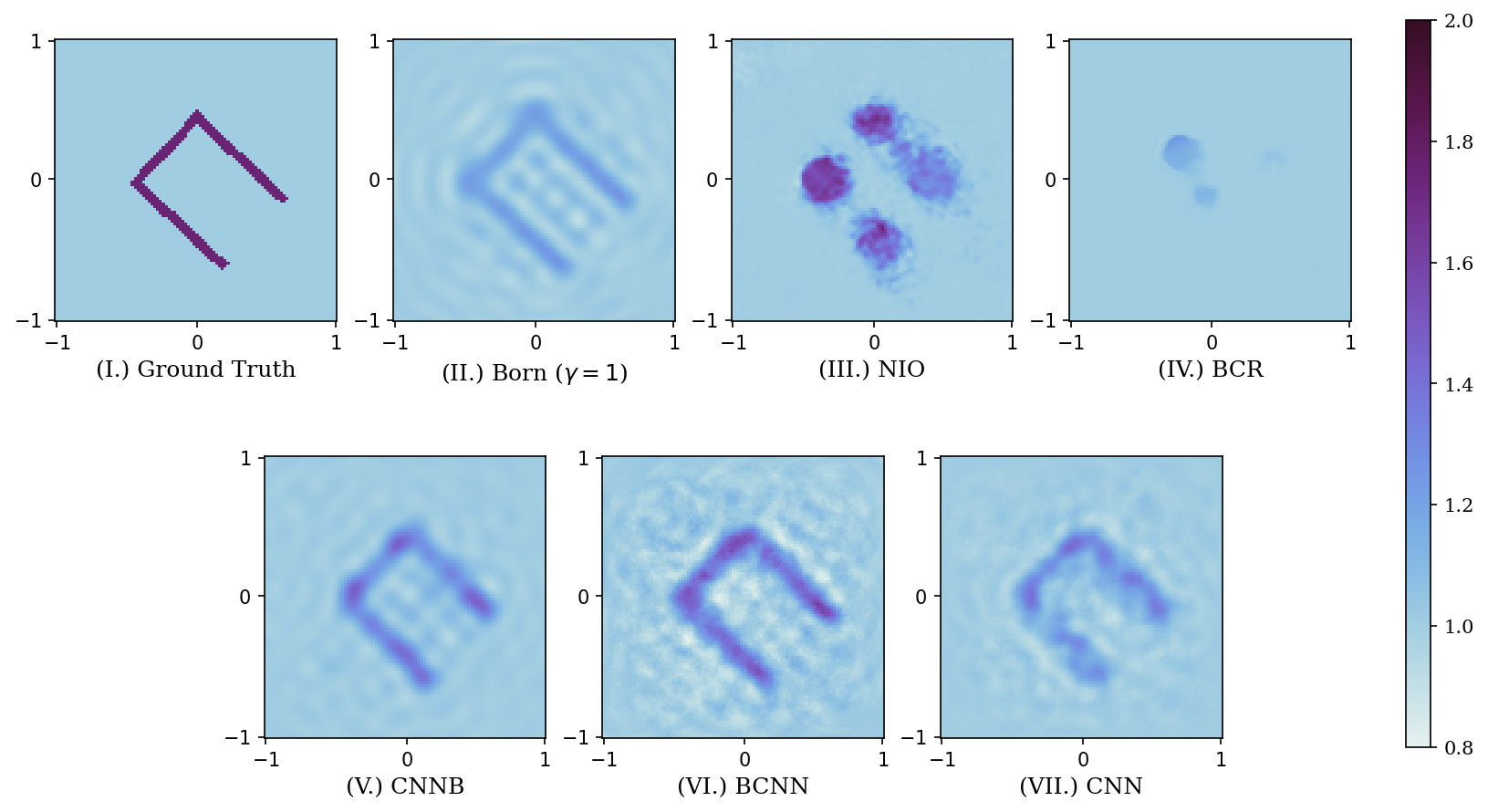}
    \caption{The reconstruction of a high-contrast, U-shaped resonance structure by the inverse Born approximation and the NN based schemes.  All three CNN based schemes result in much clearer reconstructions with CNNB having less background noise than BCNN or CNN. The NIO approach fails at constructing both the interior and the walls of the structure.}
    \label{fig:u_shape}
\end{figure}

Multiple scattering induced by resonance structures can often produce challenging far field patterns that prevent accurate inversion or reconstruction with the Born approximation. We test this using a U-shaped scatterer chosen to resonate with the incident field, as shown in Fig.\ref{fig:u_shape}. Here, we fix $\eta=1.8$ (strong scattering). 

To generate random test data,  the location, size, and rotation of each sample structure are chosen randomly from the uniform distributions $\mathcal{U}(-0.9,0.9)$, $\mathcal{U}(0.05,0.2)$, and $\mathcal{U}(0,2\pi)$, respectively (but keeping the U within $[-1,1]^2$).   In this case, results in Table~\ref{tab:errors} show that CNNB
outperforms the other methods.
To understand this, Fig.~\ref{fig:u_shape}
shows a particular realization of the inversion.  In this case,   the inverse Born method  does not approximate the the magnitude of $\eta$ well. While all three CNN models  improve on the Born approximation, the CNNB model preserves the shape of the U with minimal background speckle. Meanwhile, BCNN constructs the scatterer with distortion in the surrounding media and the CNN model is blurred. NIO produces blurred artifacts that do not resemble the U-shaped structure, while BCR-Net fails to find the scatterer.

\subsubsection{Ring}
\begin{figure}[t]
    \centering
    \includegraphics[width=0.7\linewidth]{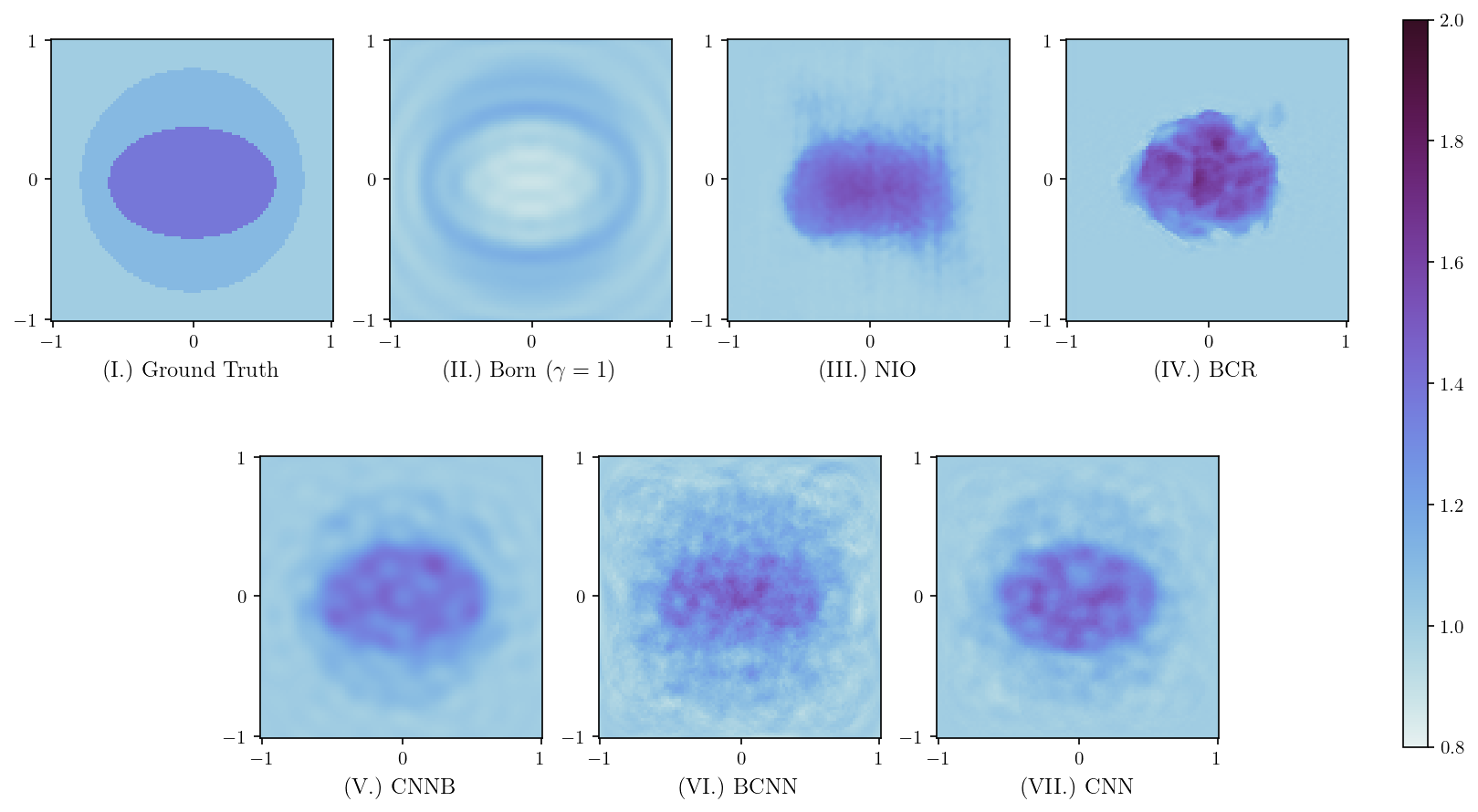}
    \caption{ The reconstruction of a high-contrast annular scatterer.  All three CNN based schemes are markedly better than inverse Born alone at representing the interior ellipses. However, the outer annular structure is poorly estimated. The NIO and BCR-Net approaches fail to detect the outer annulus.}
    \label{fig:bullseye}
\end{figure}

The goal with this example is to test if the NN based schemes can detect an inner wall, even though they were trained on scattering by simple disks. We construct an annulus or ring inside a circle with an elliptic inclusion, as shown in Fig.~ \ref{fig:bullseye}. For each sample, the value of $\eta$ in each component are chosen randomly from the uniform distribution $\mathcal{U}(1.1,1.8)$. Averaged results are shown in Table~\ref{tab:errors} showing CNNB outperforms the other methods.

For a particular choice of the model parameters, Fig.~\ref{fig:bullseye} shows that the direct use of the inverse Born scheme may not reveal the inner ellipse. In this case all three CNN models improve the inverse Born approximation and perform relatively similarly in obtaining the shape. Meanwhile, NIO and BCR-Net do not produce an accurate representation of the scatterer.

\subsubsection{Shepp-Logan Phantom}
\begin{figure}[h]
    \centering
    \includegraphics[width=0.7\linewidth]{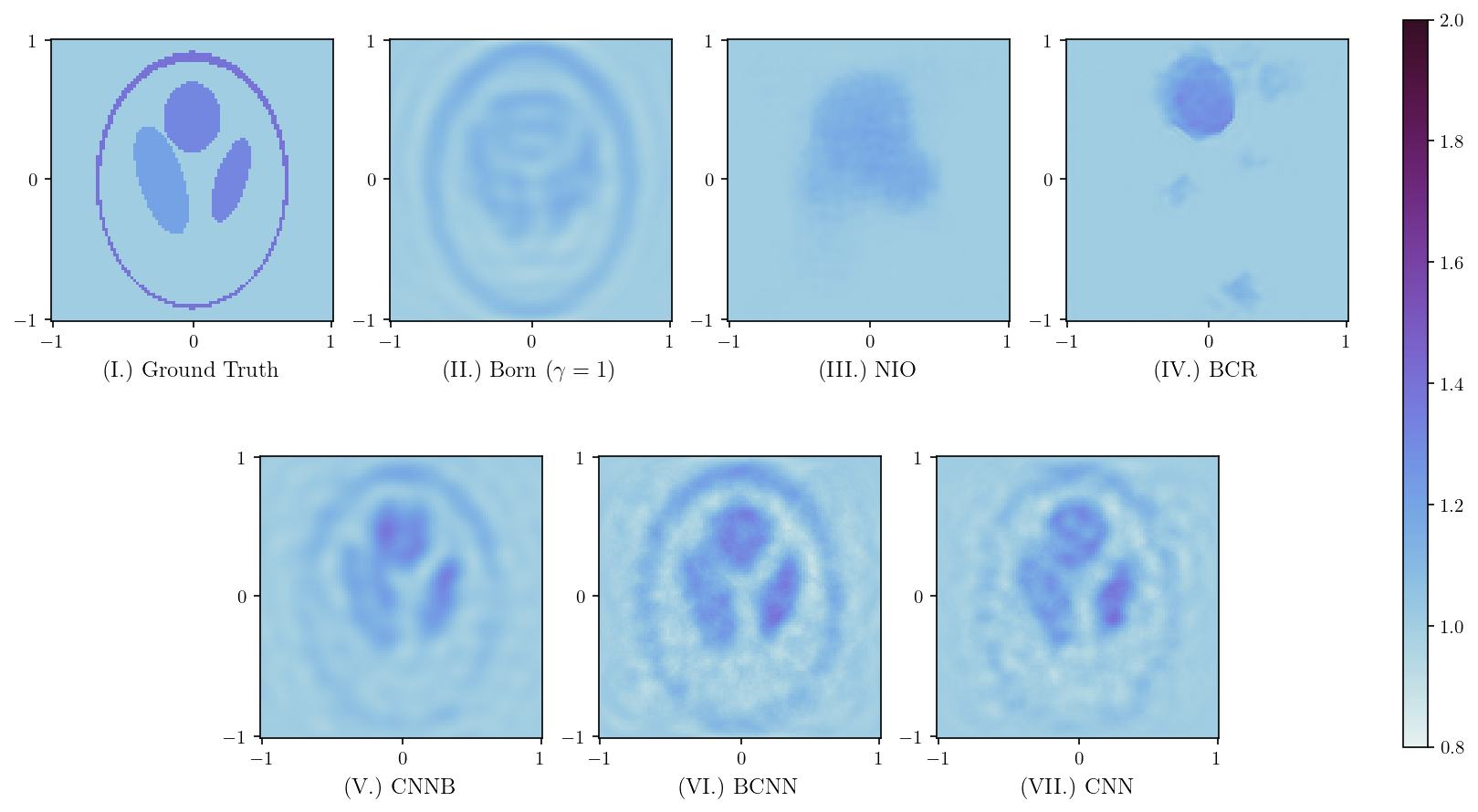}
    \caption{The reconstruction of a simplified Shepp-Logan-like phantom. Here the scattering ellipses are surrounded by a thin high contrast annulus of varying thickness.  The Born approximation can detect this ring, but fails to reconstruct inside.  All three CNN based schemes detect the three internal scatterers, though CNNB produces the sharpest image. In contrast, NIO and BCR-Net do not produce any meaningful representation of the ground truth. BCR-Net also does not capture the overall structure at all.}
    \label{fig:shepp}
\end{figure}

In potential biomedical applications, it is desirable to determine the refractive index of structures within other structures. In this example, we construct a scatterer akin to the Shepp-Logan phantom~\cite{Shepp-Logan}, as shown in Fig.~\ref{fig:shepp}. Restrictions on our mesh generator for the forward problem prevented us from using the Shepp-Logan phantom itself. To generate statistical data, for each sample the value of $\eta$ in each component is chosen randomly from the uniform distribution $\mathcal{U}(1.1,1.8)$ so that the values of $\eta$ used here are not the same as for the real Shepp-Logan phantom. Table~\ref{tab:errors} summarizes the statistical results showing that CNNB is best in the $\ell^2$ norm but CNN outperforms it in the $\ell^1$ norm.

In the particular case shown in Fig.~\ref{fig:shepp}, we see that, with appropriate regularization, the inverse Born scheme correctly images (though under-approximates) $\eta$ in the outer boundary of the scatterer, but fails to image structures inside.  All CNN based methods improve over the inverse Born approximation, though the BCNN model well-approximates both the outer ring and the separation of the three internal scatterers, CNNB gives a good quantitative reconstruction. Again, we observe that NIO and BCR-Net fail to produce an accurate representation of the phantom.

\subsubsection{Rectangles}
\begin{figure}[H]
    \centering
    \includegraphics[width=0.7\linewidth]{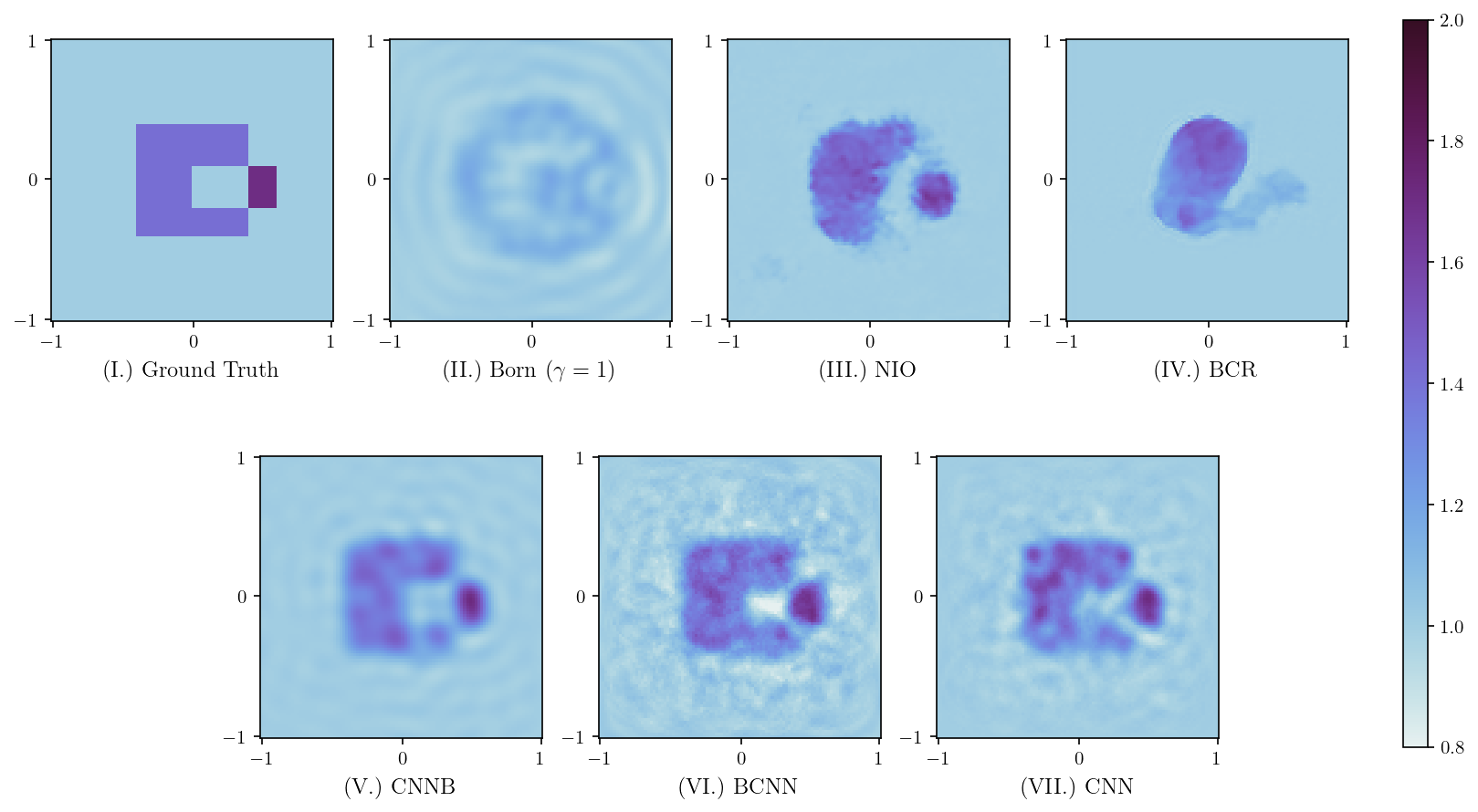}
    \caption{ The reconstruction of a structure with sharp corners.  The inner void is not seen in the Born reconstructions but is detected by the CNN schemes with CNNB and BCNN producing better reconstructions than CNN alone. NIO produces a rough inner void but  fails to handle the sharp corners of the scatterer walls.}
    \label{fig:squares}
\end{figure}

The NN based methods are trained on smooth scatterers. Our last test uses a scatterer made of rectangles and an inner void of air, as shown in Fig.~\ref{fig:squares}. For each testing sample, the value of $\eta$ in the left and right pieces are chosen randomly from the uniform distribution $\mathcal{U}(1.1,1.8)$. Table~\ref{tab:errors} shows that CNNB is, on average, the best scheme.

For the particular realization shown in
in Fig.~\ref{fig:squares} we see that the inner inclusion is invisible to the Born approximation, and the high contrast block is not well approximated. All three CNN models perform comparatively similarly, but the CNNB and BCNN models better approximate the shape of the internal region. While NIO is able to detect the interior void, it performs poorly in approximating the walls of the scatterer, while BCR-Net does not detect the void. In general, the presence of corners in this case does not badly impact the CNN based models.

\subsubsection{Discussion}
Overall, we have demonstrated good generalization to a variety of scatterers not obviously connected to the training data. These results hint that training on simple shapes like circles is a successful strategy. Moreover, we observe the quantitative success of CNNB in extending to out-of-distribution scatterers through its optimality in   $\ell^2$ and in $\ell^1$ (except for Shepp-Logan), as displayed in Table \ref{tab:errors}.

While we do not have a definitive explanation for why the CNN based schemes perform well in out-of-distribution scenarios, we conjecture that the local properties of CNNs, through their convolution blocks, may promote added structure~\cite{Goodfellow-et-al-2016}. In addition,       
the inclusion of the Born approximation in CNNB and BCNN introduces useful physical information that biases the learning process toward physically consistent reconstructions. This may improve generalization by implicitly regularizing the solution space, even beyond the weak-scattering regime.

\section{Conclusion}\label{sec:concl}

Our tests demonstrate that combining CNNs with the Born approximation can extend the applicability of the
Born approximation beyond the weak scattering limit.  Moreover, a CNN combined with the inverse Born model always improves scatterer reconstruction in comparison to a pure CNN model. 

Comparing to state of the art networks, CNNB  performs similarly to NIO and BCR-Net on in-distribution statistical tests, although we have noted that in some  cases NIO and BCR-Net fail to detect low contrast and small scatterers in close proximity to larger and higher-contrast scatterers. On noisy data, BCR-Net outperforms all the other networks although CNNB is also stable to low noise but BCNN is less stable. Both NIO and BCR-Net sometimes hallucinate at high noise levels.

When inverting data from absorbing scatterers,
NIO outperforms BCR-Net for low absorption.  It is noticeable that, in some realizations of the test data, BCR-Net incorrectly predicts 
an additional scatterer.

The CNN and Born based inversion methods excel at reconstructing the more complex out of distribution scatterers.   In this case CNNB is always best in the $\ell^2$ norm, and because of this and its stability to noise, we prefer the CNNB approach compared to CNN or  BCNN. These findings highlight that incorporating physics-based structure such as the Born approximation into learning models can improve both reconstruction accuracy and the model's ability to generalize to different data regimes. 

Much more work needs to be done to extend this
demonstration to a realistic biomedical problem like ultrasound tomography.    Fast methods should be used to compute the inverse Born approximation, such as the specialized NN of \cite{Zhou:24} or the low-rank method of \cite{Meng25}. Further, other inversion algorithms, such as an $L^1$ (LASSO) or total-variation based regularization of the Born step, could produce sharper results from CNNB at the cost of longer computing times. In addition, for specific applications, other simple training data could be considered (for example, if the goal is to image a network of blood vessels).

Another particularly interesting direction for further work is to apply the method to more exotic measurement scenarios.  For example, the case where both transmitters and receivers are on one side of the object, or when there is missing data from certain angular sectors. Moreover, additional work could include extending the applicability of the model to cases in which the wavenumber $k$ varies  (i.e. multi-frequency data) or when the scattering medium is defined by a function that is not piecewise-constant.

\section*{Acknowledgments}
 This research was supported in part through the use of the DARWIN computing system: DARWIN – A Resource for Computational and Data-intensive Research at the University of Delaware and in the Delaware Region, which is supported by NSF under Grant Number: 1919839, Rudolf Eigenmann, Benjamin E. Bagozzi, Arthi Jayaraman, William Totten, and Cathy H. Wu, University of Delaware, 2021, URL:\href{https://udspace.udel.edu/handle/19716/29071}{https://udspace.udel.edu/handle/19716/29071}. The authors acknowledge the CSC – IT Center for Science, Finland, for generously sharing their computational resources. We thank Dr. Yuwei Fan and Professor Lexing Ying for making the source code for BCR-Net available to us. We also thank Professor Ke Chen of the University of Delaware for sharing his computational resources.  Finally, we thank the anonymous referees for their very valuable suggestions.

\bibliographystyle{siamplain}
\bibliography{NN-born}

\end{document}